\documentclass[12pt]{amsart}

\usepackage{latexsym,amsmath,a4wide,amssymb,amsthm,graphics}

\usepackage{comment}
\usepackage{nicematrix}
\usepackage{multirow}

\usepackage{tikz,rotating}
\usetikzlibrary{matrix}
\usepackage{enumitem}
\usepackage{arydshln}
\usepackage{setspace}
\usepackage{ifpdf}

\usepackage{ytableau}
\usetikzlibrary{tikzmark}
\usepackage{nicematrix}

\usepackage{tikz-cd}

\usepackage{amsmath,lipsum}

\usepackage[utf8]{inputenc}
\usepackage[OT2,T2A]{fontenc}
\usepackage{graphicx}
\usepackage{gensymb}
\usepackage{enumitem}
\usepackage{setspace}
\usepackage{tabularx}
\usepackage{upgreek}
\usepackage{hyperref}
\hypersetup{colorlinks=true, urlcolor=cyan, linkcolor=black}
\newcommand{\ran}{\text{ran}}

\newcommand{\CC}{\mathbb{C}}
\newcommand{\NN}{\mathbb{N}}
\newcommand{\ZZ}{\mathbb{Z}}

\newcommand{\cD}{\overline{\mathbb{D}}}

\theoremstyle{definition}
\newtheorem{remark}{Remark}[section]

\newtheorem*{thm}{Theorem 5.1}
\newtheorem{theorem}{Theorem}[section]
\newtheorem{definition}{Definition}[section]
\newtheorem{proposition}{Proposition}[section]

\newtheorem{corollary}{Corollarly}[section]

\title{The Taylor spectrum of pairs of isometries}
\author{Zbigniew Burdak}
\author{Patryk Pagacz}
\address[Zbigniew Burdak]{Department of Applied Mathematics\\
University of Agriculture in Krakow\\
Balicka 253C\\
30-198 Krak\'ow\\
Poland}

\email[Zbigniew Burdak]{zburdak@urk.edu.pl\\rmburdak@cyf-kr.edu.pl}

\address[Patryk Pagacz]{Faculty of Mathematics and Computer Science\\
   Jagiellonian University\\
   \L ojasiewicza 6\\
   30-348 Krak\'ow\\
   Poland}
\email[Patryk Pagacz]{patryk.pagacz@uj.edu.pl\\patryk.pagacz@gmail.com }
\begin{document}
\maketitle

\begin{abstract}
In the paper we fully describe Taylor spectrum of pairs of isometries given by diagrams. In most cases both isometries in such pairs have non-trivial shift part and its Taylor spectrum is a proper subset (of Lebesgue measure in $(0,\pi^2)$) of the closed bidisc. 
\end{abstract}
\subjclass{MSC2020: 47B02, 47B40, 47A13}

\keywords{Keywords: a pair of isometries, spectrum, Taylor spectrum}

\section{Introduction}

Spectrum of an operator and functional calculus are fundamental concepts in Functional Analysis. Apparently simple examples of a unilateral shift or a bilateral shift are enough to see how much more complex is the concept of the spectrum of an operator comparing to the final dimensional case - the set of eigenvalues (spectrum) of a matrix. In many particular cases the spectrum of an operator can be determined. For example the spectrum of an isometry can be a closed unit disc or a closed subset of a unit circle (a set of Lebesgue measure zero). Much less is known about the spectrum of pairs/tuples of commuting operators.
Ideas of generalization of the concept of spectrum to  pairs/tuples of commuting operators have been appearing since 50's (see \cite{Arens,CurtoSurvey,Dash1,Dash2,Harte,Slodkowski,Taylor,Zelazko}). The best one seems to be the definition of Taylor spectrum. It is mostly because there exist functional calculus for
functions analytic on a neighbourhood of the Taylor spectrum (see e.g. \cite{Albrecht,CurtoSurvey,MulCalc,TaylorCalc,VCalc}). For this reason
many experts consider the Taylor spectrum to be the proper generalization of the
 single operator spectrum (see \cite{Muller},  \cite[Table 2.6]{CurtoSurvey}).
Additionally, Taylor spectrum technique finds their application in studies e.g. on specific classes of subnormal operators, see \cite{Stochel}.

Unfortunately Taylor spectrum is not easy to be calculated. The best example is still missing Taylor spectrum of pairs of commuting isometries, which partial answer is the subject of this paper.
There are some tools (see e.g. \cite{CCLJ,Curto0,V}) which however have limitations. Let us give (not exhausting) list of contributions which have been made over the years.
The spectrum of left and right multiplications induced by some tuple of operators was consider in \cite{Curto2}.
Next, Curto and Yan (see \cite{CY1,CY2}) determinate the Taylor spectrum of the multiplications by the coordinates functions belonging to  the closure of polynomials in $L^2(\mu)$, for Reinhardt measure $\mu$.
Such operators are commuting weighted shifts, i.e. the pairs $(T_1,T_2)$ such that $T_ie_{(k,l)}=\omega_{k,l}e_{(k,l)+\varepsilon_i}$, for $k,l\geq0$ and some $\omega_{k,l}\in\CC$, where $\varepsilon_1=(1,0)$ and $\varepsilon_2=(0,1)$.
Another class of weighted shifts was considered in \cite{Curto}.
Namely, the Taylor spectrum of hyponormal 2-variable weighted shifts with commuting subnormal components was determined there.

In the recent paper \cite{Bara} Bhattacharyya et al. considered Taylor spectrum of pairs of commuting isometries. They calculated the spectrum for the cases where defect operator (see \cite{Guo, He}) is non-negative (which is doubly commuting case) or non-positive (the only pair of non-unitary operators with non-positive defect is the modified bi-shift - see \cite{Pop1,Pop2}). The non-positive case was later used in \cite{Bera} to consider $2$-isometries.
Recall that by \cite{Slo}, the only doubly commuting shifts are multiplications by coordinates $M_w,M_z$ on the Hardy space.  The both classes, doubly commuting shifts and modified bi-shifts are certain cases of pairs of isometries defined by a diagram introduced in \cite{Mil}. In this paper we calculate the spectrum of a general pair of isometries given by a diagram.
The class is a crucial type of compatible isometries (see \cite{Comp,Comp2}) and in most cases their defect operator is a difference of mutually orthogonal projection. Moreover, it
covers the class of totally non-unitary pair of isometries, introduced in \cite{S}. Therefore, it plays a crucial role in the study of the weak-stationary stochastic processes (see \cite{procesy}).

The main result of the paper is 
\begin{thm}
  Let $M_w, M_z\in \mathcal{B}(H_J)$, where $J$ is a non-simple diagram.

    Then $$\sigma_T(M_w,M_z)= \{(\mu,\lambda)\in \overline{\mathbb{D}}^2:
|\mu|^{\max(\rho_-,\rho_+)}\leq|\lambda|\leq |\mu|^{\min(\delta_-,\delta_+)}\},$$ 
where $\delta_-,\delta_+,\rho_-,\rho_+$ are defined by \eqref{deltarho+-} and inequalities with $0^0$ or $1^\infty$ are assumed to be satisfied.

\end{thm}
where $M_w, M_z\in \mathcal{B}(H_J)$ is a model of a pair of isometries given by a diagram, and non-simple diagram is such that the corresponding pair is not doubly commuting. The parameters $\delta_-,\delta_+,\rho_-,\rho_+$ are in $[0,\infty]$ and depends on the shape of the diagram only. Moreover, for any $0\leq a\leq b\leq\infty$ there is a diagram such that $\min(\delta_-,\delta_+)=a, \max(\rho_-,\rho_+)=b$.

\section{Preliminaries}
Let $\NN=\{0,1,\ldots\},\ZZ_+=\{1,2,\ldots\},\ZZ_-=\{-1,-2,\ldots\}$, $\mathbb{D}=\{z\in\mathbb{C}:|z|< 1\}$ with closure $\cD=\{z\in\mathbb{C}:|z|\leq 1\}$ and border $\partial
\mathbb{D}=\mathbb{T}=\{z\in\mathbb{C}:|z|= 1\}.$ Further, $B(\mathcal{H})$ denotes the algebra of bounded, linear operators on a Hilbert space $\mathcal{H}.$ The standard notation $L^2(\mathbb{T}), L^2_{\mathcal{H}}(\mathbb{T}),L^2(\mathbb{T}^2), L^2_{\mathcal{H}}(\mathbb{T}^2)$ is for the spaces of square integrable, scalar or $\mathcal{H}$ valued functions over circle $\mathbb{T}$ or torus $\mathbb{T}^2,$ respectively.  Operators of multiplication for the scalar valued case are denoted $M_z$ over circle and $M_w, M_z$ over torus. Since $L^2_{\mathcal{H}}(\mathbb{T})$ ($L^2_{\mathcal{H}}(\mathbb{T}^2)$) is unitarily equivalent with $L^2(\mathbb{T})\otimes\mathcal{H}$ ($L^2(\mathbb{T}^2)\otimes\mathcal{H}$) operators of multiplication in vector valued case are denoted $M_z\otimes I_{\mathcal{H}}$ and $M_w\otimes I_{\mathcal{H}}, M_z\otimes I_{\mathcal{H}}$ respectively. The space of square summable sequences valued in $\mathcal{H}$ is denoted by $l^2(\mathcal{H})$ and $l^2$ if $\mathcal{H}=\mathbb{C}.$

The following remark justifies the convention to interpret inequalities with indeterminate form $0^0$ or $1^\infty$ as satisfied. 

\begin{remark}\label{pab}
    Let $0\leq p\leq q\leq \infty$ and 
    \begin{equation}\tag{$\dagger$}\label{one}a^q\leq b\leq a^p \text{ for } a\in(0,1), b\in[0,1], \quad \text{ or }\quad b^\frac{1}{p}\leq a \leq b^\frac{1}{q}\text{ for } a\in[0,1], b\in (0,1),\end{equation} where we assume $ \frac{1}{0}=\infty$ and as usual take  $\frac{1}{\infty}=0.$

    On the other hand, let 
    \begin{equation}\tag{$\ddagger$}\label{two}a^q\leq b\leq a^p \text{ for } a,b\in[0,1]\end{equation}
    where inequalities with the symbol $0^0$ or $1^\infty$ are assumed to be satisfied. 
    
    One can check that \eqref{one} is equivalent to \eqref{two} on $(0,1)\times[0,1]\cup[0,1]\times(0,1).$ Moreover, pairs $(0,0), (1,1)$ satisfy \eqref{two} for all $p,q,$ while $(0,1),(1,0)$ satisfy \eqref{two} only for $p=0, q=\infty$, respectively.
    
    Let us summarize the condition \eqref{two} in all the cases where at least one of $p, q$ is $0$ or $\infty$.

The condition \eqref{two} is equivalent to:    
     \begin{itemize}
          \item $a^q\leq b$ for $0=p<q<\infty,$
        \item $b\leq a^p$ for $0<p<q=\infty,$
         \item $a=0$ or $b=1$ for $0=p=q,$
         \item $a=1$ or $b=0$ for $p=q=\infty,$ 
         \item $[0,1]^2$ for $p=0, q=\infty.$ 
    \end{itemize}
    \end{remark}
    \begin{proof}
 For $0<p\leq q<\infty$ both pairs of inequalities in \eqref{one} make sense on $[0,1]^2$ and such extensions are equivalent to each other, so to \eqref{two}. The cases where at least one of $p, q$ is $0$ or $\infty$  are listed in the summarize part. The equivalence of \eqref{one} and \eqref{two} on $(0,1)\times[0,1]\cup[0,1]\times(0,1)$ in each of the cases and the remaining part of remark may be checked by a direct calculation.
 \end{proof}
\subsection{Spectrum of a single operator}
We use the standard notation $\sigma(T)$ for the spectrum of a single operator $T\in B(H)$ and $\sigma_p(T), \sigma_c(T), \sigma_r(T), \sigma_{ap}(T)$ for point spectrum, residual spectrum, approximate point spectrum and continuous spectrum, respectively. Namely,
\begin{itemize}
    \item $\sigma(T)=\{\lambda\in \CC: \lambda- T \textnormal{ is not invertible} \},$
    \item $\sigma_p(T)=\{\lambda\in \CC: \ker(\lambda- T)\not= \{0\} \},$
    \item $\sigma_c(T)=\{\lambda\in \CC\setminus\sigma_p(T): \ran(\lambda- T)\not=\mathcal{H}, \overline{\ran(\lambda- T)}=\mathcal{H}\},$
    \item $\sigma_r(T)=\{\lambda\in \CC\setminus\sigma_p(T): \overline{\ran(\lambda- T)}\neq \mathcal{H}\}=\sigma(T)\setminus (\sigma_p(T)\cup\sigma_c(T))$,
    \item $\sigma_{ap}(T)=\{\lambda\in\CC: \|(\lambda-T)x_n\|\to 0 \textnormal{ for some } x_n \textnormal{ such that } \|x_n\|=1\}$.
\end{itemize}

Let us recall some well known results connected with the spectrum of a single operator.
\begin{remark}\label{Clcond}
   The range of the operator $A\in B(\mathcal{H},\mathcal{K})$ is closed if and only if there is $c>0$ such that $$\|Ax\|\geq c\|x\|,$$ for all $x\bot \ker A$.
\end{remark}
We take advantage of some properties of spectrum of a unilateral shift.
\begin{remark}\label{RS_is_cl}
    Let $S\in B(l^2)$ be a unilateral shift. Then:
    \begin{itemize}
        \item $\sigma_p(S)=\emptyset, \sigma_r(S)=\mathbb{D}, \sigma_{ap}(S)=\sigma_c(S)=\mathbb{T}.$
        \item $\ran(\lambda-S)$ is closed for $\lambda\in\mathbb{D}$ and $\ker (\overline{\lambda}-S^*)=\mathbb{C}\vec{v}$ where $\vec{v}=(1,\overline{\lambda},\overline{\lambda}^2,\dots).$
        \item $\ran(\lambda-S)$ is not closed for $\lambda\in\mathbb{T}.$
    \end{itemize}
\end{remark}

\subsection{Taylor spectrum of a pair of operators}
Let us consider a commuting pair of operators $\textbf{T}=(T_1,T_2)$ 
and a short Koszul complex $K(\textbf{T}, \mathcal{H})$ associated to $\textbf{T}$ on $\mathcal{H}$:
$$K(\textbf{T}, \mathcal{H}) : \{0\} \xrightarrow{0_1} \mathcal{H}\xrightarrow{\delta_1}\mathcal{H}\oplus \mathcal{H}\xrightarrow{\delta_2}\mathcal{H}\xrightarrow{0_2} \{0\},$$
where $\delta_1(h)=(T_1h,T_2h)$ and $\delta_2(h_1,h_2)=-T_2h_1+T_1h_2$. Since $\delta_2 \circ \delta_1 =0$ the Koszul complex is well defined.
We say that the Koszul complex is \emph{exact} 
if $\ker \delta_1 = \ran (0_1)$, $\ker \delta_2=\ran (\delta_1)$ and $\ran (\delta_2) = \ker 0_2$.
In other words, we say that $K(\textbf{T}, \mathcal{H})$ is exact if
\begin{enumerate}[label=(T\arabic*)]
    \item\label{T1} $\ker T_1 \cap \ker T_2 = \{0\}$, 
    \item\label{T2} for any $h_1,h_2\in \mathcal{H}$ such that $T_1h_2-T_2h_1 =0$ there is $h\in \mathcal{H}$ such that $$h_1=T_1h,\quad \textnormal{ and } h_2=T_2h,$$
    \item\label{T3} $\ran(T_1) + \ran(T_2) = \mathcal{H}$.
\end{enumerate}

The exactness of Koszul complex is much less then invertibility of both operators. In fact invertibility of one of operators is enough. 
\begin{remark}\label{inver_koszul}
   If any of $T_1, T_2$ is invertible, then Koszul complex is exact. 
\end{remark}
\begin{proof}
    Conditions \ref{T1}, \ref{T3} are immediate. Let $T_1$ be the one invertible and $T_1h_2-T_2h_1=0.$ Obviously \ref{T2} may be satisfied only with $h=T_1^{-1}h_1.$ It is indeed, as $T_1h_2-T_2h_1=0$ implies $T_1h_2=T_2h_1$ and in turn $h_2=T_1^{-1}T_2h_1=T_2T_1^{-1}h_1=T_2h.$
\end{proof}
\begin{definition}
\emph{Taylor spectrum} of a commuting pair $(T_1,T_2)$ of operators on the Hilbert space $\mathcal{H}$ is the set
$$\sigma_T(T_1,T_2)=\{(\lambda_1,\lambda_2)\in \CC^2 : K((\lambda_1-T_1,\lambda_2-T_2),\mathcal{H}) \textnormal{ is not exact} \}.$$\end{definition}
It is known that Taylor spectrum is closed (see f.e. \cite[chap.IV, Theorem 4]{Muller}) and \begin{equation}\label{tsinsts}\sigma_T(T_1, T_2)\subset\sigma(T_1)\times\sigma(T_2).\end{equation} The latter inclusion follows from Remark \ref{inver_koszul}. Since exactness of Koszul complex may brake at any of stages \ref{T1}, \ref{T2}, \ref{T3} we get a natural decomposition $\sigma_T(T_1,T_2)=\Gamma_1(T_1,T_2)\cup \Gamma_2(T_1,T_2)\cup \Gamma_3(T_1,T_2)$ where $$\Gamma_{\textnormal{i}}(T_1,T_2):=\{(\lambda_1,\lambda_2): K((\lambda_1-T_1,\lambda_2-T_2),\mathcal{H}) \textnormal{ does not satisfy (Ti)} \}.$$
Similar idea of investigating Taylor spectrum in such parts  appeared in 
\cite{Muller} (where completions of such sets were investigated) and \cite{Bara}.
\begin{remark}\label{T2iso}
Let $T_1, T_2\in B(\mathcal{H})$ and $\lambda_1\notin \sigma_p(T_1)$. Then the pair $(\lambda_1-T_1, \lambda_2-T_2)$:
\begin{itemize}
    \item satisfy \ref{T1}, 
    \item condition \ref{T2} is equivalent to:
         \begin{equation}\tag{T2'}\label{T2'}(\lambda_2-T_2)h_1\in\ran(\lambda_1-T_1)\Rightarrow h_1\in \ran{(\lambda_1-T_1)}, \text{ for all } h_1\in \mathcal{H}.\end{equation}
\end{itemize}
\end{remark}
\begin{proof}
The pair satisfy  \ref{T1} as $\ker (\lambda_1-T_1)=\{0\}$ by $\lambda_1\notin \sigma_p(T_1)$.

For the second part note that $(\lambda_2-T_2)h_1\in\ran(\lambda_1-T_1)$ is equivalent to $(\lambda_1-T_1)h_2-(\lambda_2-T_2)h_1 =0$ for some $h_2,$ so \eqref{T2'} is equivalent to:
 $$\text{for any } h_1,h_2\in \mathcal{H}\text{ such that }(\lambda_1-T_1)h_2-(\lambda_2-T_2)h_1 =0\text{ the vector }h_1\in \ran{(\lambda_1-T_1)}.$$
 
 Thus \ref{T2} implies \eqref{T2'}. For the reverse implication assume \eqref{T2'} and take $h_1,h_2\in \mathcal{H}$ such that $(\lambda_1-T_1)h_2-(\lambda_2-T_2)h_1 =0,$ so by \eqref{T2'} we get $h_1\in \ran{(\lambda_1-T_1)}.$ Let $h_1=(\lambda_1-T_1)h$ for some vector $h$. Then $(\lambda_1-T_1)h_2=(\lambda_2-T_2)h_1$ implies $$(\lambda_1-T_1)(h_2-(\lambda_2-T_2)h)=(\lambda_2-T_2)h_1-(\lambda_2-T_2)(\lambda_1-T_1)h=(\lambda_2-T_2)(h_1-(\lambda_1-T_1)h)=0.$$ Since $\lambda_1\notin \sigma_p(T_1)$ we get $h_2=(\lambda_2-T_2)h.$
\end{proof}

Note, that $\Gamma_1(T_1,T_2)\subset\sigma_p(T_1)\times\sigma_p(T_2)$ and the inclusion is strong in general. Indeed, \ref{T1} fails only if eigenspaces corresponding to the pair of eigenvalues has nontrivial intersection.

\begin{remark}\label{cond3}
The condition \ref{T3} is not satisfied in two (not necessarily disjoint) cases:
\begin{itemize}
    \item $\ker  T^*_1 \cap \ker T^*_2 \not= \{0\},$
    \item $\ran(T_1) + \ran(T_2)$ is not closed.
\end{itemize}

Indeed, the condition $\overline{\ran(T_1) + \ran(T_2)} \not= \mathcal{H}$ means that there is $h\in \mathcal{H}$ such that for any $x,y\in \mathcal{H}$ we get
$\langle T_1x+T_2y,h\rangle=0$.
So $\langle x,T_1^*h\rangle+\langle y,T_2^*h\rangle=0$. In other words, $T_1^*h=T_2^*h=0$.
\end{remark}

\subsection{Pairs of isometries defined by a diagram}
We recall the concept of pairs of isometries defined by a diagram. 
Let us start with definition of a diagram with corresponding notions. 
\begin{definition}\cite{Mil, Comp, BKPS}\label{diagram definition}
 A set $J \subset \mathbb{Z}^2$ is called a \emph{diagram} if $J+\mathbb{N}^2\subset J$, that is $(i+k, j+l)\in J$ for any $(i,j)\in J, (k,l) \in \mathbb{N}^2$ .

 Diagrams $J_1, J_2$ are (translation) equivalent if $J_1=(i,j)+J_2$ for some $(i,j)\in\mathbb{Z}^2.$

 Diagram equivalent to any of $\mathbb{Z}^2, \mathbb{Z}_+^2, \mathbb{Z}_+\times \mathbb{Z}, \mathbb{Z}\times \mathbb{Z}_+$ is called simple.

 The sets $J_{VB} = \{(i,j) \in J : (i-1,j) \not\in J \}$ and
$J_{HB} = \{(i,j) \in J : (i,j -1) \not\in J \}$ are called the \emph{vertical} and \emph{horizontal borders} of $J.$

The set of outer corners of a diagram is $J_{outer}=J_{HB}\cap J_{VB}.$

The set of inner corners is $J_{inner}=\{(i,j)\in J: (i-1,j)\in J, (i,j-1)\in J, (i-1, j-1)\notin J\}$.

Border of the diagram is the set $J_{border}=J_{HB}\cup J_{VB}\cup J_{inner}$
\end{definition}

For $j\in \mathbb{Z}$ and a diagram $J$ define $M_j=\infty$ if $\{i\in\mathbb{Z}:(i,j)\in J\}=\emptyset$ and $M_j:=\inf\{i\in\mathbb{Z}:(i,j)\in J\}$ otherwise.   Note that $J+\mathbb{N}^2\subset J$ and definition of $M_j$  yields $J\cap (\mathbb{Z}\times\{j\})=\{(i,j)\in\mathbb{Z}^2:M_j\leq i\}$ for any $j.$ In particular, $M_j=-\infty$ if $\mathbb{Z}\times \{j\}\subset J.$ Moreover, \begin{equation}\label{Jform}
J=\bigcup_{j=-\infty}^\infty\{(i,j)\in \ZZ^2: M_j\leq i\}.    
\end{equation} 
Since $J+\mathbb{N}^2\subset J$ yields $(M_{j-1},j-1)+(0,1)=(M_{j-1},j)\in J$ we get $M_j\leq M_{j-1}$ for any $j\in\mathbb{Z}$. In other words, $\{M_j\}_{j=-\infty}^\infty$ is a non-increasing sequence in $\ZZ\cup\{-\infty, \infty\}$. Note also, that
\begin{equation}\label{d1}J_{VB}=\{(M_j,j): j\in\mathbb{Z}\text{ such that }M_j\in\mathbb{Z}\}\end{equation}
and 
\begin{equation}\label{d2}J_{HB}=\bigcup_{j:M_j<M_{j-1}}\{(i,j)\in\mathbb{Z}^2:M_j\leq i<M_{j-1}\}.\end{equation}

In the same way, there is defined a non-increasing sequence $\{N_i\}_{i=-\infty}^\infty$ such that 
\begin{equation}\label{JformN}
J=\bigcup_{i=-\infty}^\infty\{(i,j)\in \ZZ^2: N_i\leq j\}.    
\end{equation}
with the corresponding description of borders:
$$J_{VB}=\bigcup_{i:N_i<N_{i-1}}\{(i,j)\in\mathbb{Z}^2:N_i\leq j<N_{i-1}\},\quad  J_{HB}=\{(i,N_i): i\in\mathbb{Z}\text{ such that } N_i\in\mathbb{Z}\}.$$
Let us picturize a diagram and the corresponding concepts.

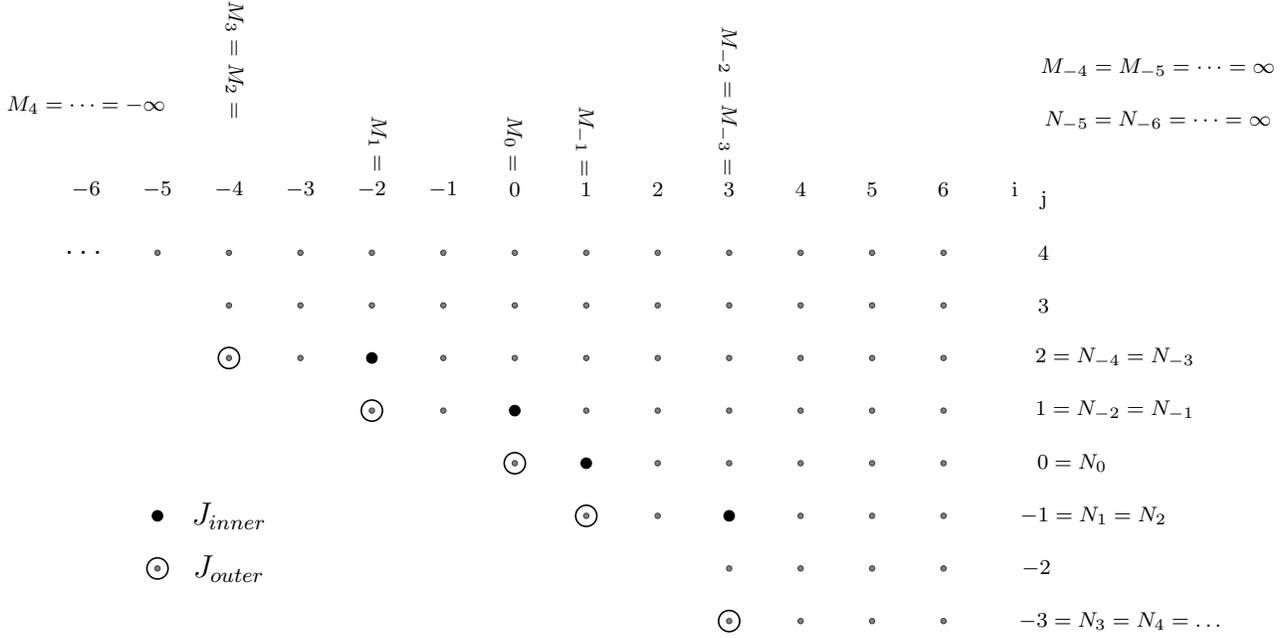
\begin{center}
 \begin{tikzpicture}
[yscale=0.7, xscale=0.95, auto,
kropka/.style={draw=black!75, circle, inner sep=0pt, minimum size=2pt, fill=black!50},
kropa/.style={draw=black, circle, inner sep=0pt, minimum size=4pt, fill=black!100},
kolo/.style={draw=black, circle, inner sep=0pt, minimum size=8pt, semithick},
kwadrat/.style={draw=black, rectangle, inner sep=0pt,  minimum size=10pt, semithick},
romb/.style={draw=black, diamond, inner sep=0pt, minimum size=8pt, semithick},
trojkat/.style ={draw=black, regular polygon, regular polygon sides=3, inner sep=0pt, minimum size=8pt, semithick}]

\path (6,-2) node [kropka] {} (5,-2) node [kropka] {} (4,-2) node [kropka] {};
\path (6,-1) node [kropka] {} (5,-1) node [kropka] {} (4,-1) node [kropka] {};
\path (6,0) node [kropka] {} (5,0) node [kropka] {} (4,0) node [kropka] {} (3,0) node [kropka] {} (2,0) node [kropka] {};
\path (6,1) node [kropka] {} (5,1) node [kropka] {} (4,1) node [kropka] {} (3,1) node [kropka] {} (2,1) node [kropka] {}
(1,1) node [kropka] {};
\path (6,2) node [kropka] {} (5,2) node [kropka] {} (4,2) node [kropka] {} (3,2) node [kropka] {} (2,2) node [kropka] {}
(1,2) node [kropka] {} (0,2) node [kropka] {} (-1,2) node [kropka] {};
\path (6,3) node [kropka] {} (5,3) node [kropka] {} (4,3) node [kropka] {} (3,3) node [kropka] {} (2,3) node [kropka] {}
(1,3) node [kropka] {} (0,3) node [kropka] {} (-1,3) node [kropka] {} (-2,3) node [kropka] {} (-3,3) node [kropka] {};
\path (6,4) node [kropka] {} (5,4) node [kropka] {} (4,4) node [kropka] {} (3,4) node [kropka] {} (2,4) node [kropka] {}
(1,4) node [kropka] {} (0,4) node [kropka] {} (-1,4) node [kropka] {} (-2,4) node [kropka] {} (-3,4) node [kropka] {} (-5,4) node [kropka] {} (-6,4) node[]{$\dots$};

\path (6,-3) node [kropka] {}  ++(-1,0) node [kropka]  {}
        ++(-1,0) node [kropka] {} ++(-1,0) node [kolo] {} ++(0,0) node [kropka] {} ++(0,1) node [kropka] {}
        ++(0,1) node [kropa] {} ++(-1,0) node [kropka] {}
        ++(-1,0) node [kolo] {} ++(0,0) node [kropka] {} ++(0,1) node [kropa] {}
        ++(-1,0) node [kolo] {} ++(0,0) node [kropka] {} ++(0,1) node [kropa] {}
        ++(-1,0) node [kropka] {}  ++(-1,0) node [kolo] {} ++(0,0) node [kropka] {}
        ++(0,1) node [kropa] {} ++(-1,0) node [kropka] {} ++(-1,0) node [kolo] {} ++(0,0) node [kropka]{}
        ++(0,1) node [kropka] {} ++(0,1) node [kropka] {};
\path (-5,-2) node [kolo] {} ++(0,0) node [kropka] {};
\node at (-4,-2){$J_{outer}$};
\path (-5,-1) node [kropa] {};
\node at (-4,-1){$J_{inner}$};

\node at (-6,6.8) {\tiny{$M_4=\dots=-\infty$}};
\node at (-6,5.2) {\tiny{$-6$}};
\node at (-5,5.2) {\tiny{$-5$}};
\node at (-4,8.8) {\begin{rotate}{-90}\tiny{$M_3=M_2=$}\end{rotate}};
\node at (-4,5.2) {\tiny{$-4$}};
\node at (-3,5.2) {\tiny{$-3$}};
\node at (-2,6.6) {\begin{rotate}{-90}\tiny{$M_1=$}\end{rotate}};
\node at (-2,5.2) {\tiny{$-2$}};
\node at (-1,5.2) {\tiny{$-1$}};
\node at (-0.1,6.6) {\begin{rotate}{-90}\tiny{$M_0=$}\end{rotate}};
\node at (0,5.2) {\tiny{$0$}};
\node at (0.9,6.8) {\begin{rotate}{-90}\tiny{$M_{-1}=$}\end{rotate}};
\node at (1,5.2) {\tiny{$1$}};
\node at (2,5.2) {\tiny{$2$}};
\node at (2.9,8.3) {\begin{rotate}{-90}\tiny{$M_{-2}=M_{-3}=$}\end{rotate}};
\node at (9,7.5) {\tiny{$M_{-4}=M_{-5}=\dots=\infty$}};
\node at (9,6.5) {\tiny{$N_{-5}=N_{-6}=\dots=\infty$}};
\node at (3,5.2) {\tiny{$3$}};
\node at (4,5.2) {\tiny{$4$}};
\node at (5,5.2) {\tiny{$5$}};
\node at (6,5.2) {\tiny{$6$}};
\node at (7,5.2) {\tiny{i}};

\node at (7.4, 5) {\tiny{j}};
\node at (7.4, 4) {\tiny{$4$}};
\node at (7.4, 3) {\tiny{$3$}};
\node at (8.4, 2) {\tiny{$2=N_{-4}=N_{-3}$}};
\node at (8.4, 1) {\tiny{$1=N_{-2}=N_{-1}$}};
\node at (7.8, 0) {\tiny{$0=N_0$}};
\node at (8.1, -1) {\tiny{$-1=N_1=N_2$}};
\node at (7.3, -2) {\tiny{$-2$}};
\node at (8.5, -3) {\tiny{$-3=N_3=N_4=\dots$}};

\end{tikzpicture}
\end{center}

\begin{definition}
\label{on_Diagram}
Let $J$ be a diagram and $\mathcal{H}$ be a Hilbert space. Define $$H_J(\mathcal{H}):=\{f\in L^2_{\mathcal{H}}(\mathbb{T}^2): \hat{f}_{i,j}=0\text{ for }(i,j)\notin J\}=\overline{\text{Span}\{w^iz^jh:(i,j)\in J, h\in\mathcal{H}\}}$$ where  $\hat{f}_{i,j}$ denotes Fourier coefficients. Note that $H_J(\mathcal{H})$ is invariant under operators of multiplication by independent variables $M_w\otimes I_{\mathcal{H}}, M_z\otimes I_{\mathcal{H}}$.

A pair of isometries defined by the diagram $J$ and the space
$\mathcal{H}$ is a pair unitarily equivalent to $(M_w\otimes I_{\mathcal{H}}|_{H_J(\mathcal{H})}, M_z\otimes I_{\mathcal{H}}|_{H_J(\mathcal{H})})=(M_w|_{H_J}\otimes I_{\mathcal{H}},M_z|_{H_J}\otimes I_{\mathcal{H}}),$ where $H_J=H_J(\mathbb{C})\subset L^2(\mathbb{T}^2).$

A pair of isometries defined by the diagram $J$ is the pair of isometries defined by the diagram $J$ and any space. 
\end{definition}
We may limit to scalar valued case by the following observation.
\begin{remark}
    The Taylor spectrum of the pair given by a diagram depends on the diagram only. More precisely $$\sigma_T(M_w|_{H_J}\otimes I_{\mathcal{H}}, M_z|_{H_J}\otimes I_{\mathcal{H}})=\sigma_T(M_w|_{H_J}, M_z|_{H_J})$$ for any $\mathcal{H}\neq\{0\}.$
\end{remark}

In the paper we denote $M_w, M_z\in B(H_J),$  instead of $M_w|_{H_J}, M_z|_{H_J}.$

The defect operator $C(V_1,V_2)$ of a pair of commuting isometries $(V_1,V_2)$ was introduced in \cite{Guo, He}
$$C(V_1,V_2)=I-V_1V_1^*-V_2V_2^*+V_1V_2V_2^*V_1^*.$$ It is clear that $C(V_1,V_2)$ vanishes on $\ran(V_1V_2)$ and $$\ker V_1^*V_2^*=\ker V_1^*\oplus V_1(\ker V_2^*)=\ker V_2^*\oplus V_2(\ker V_1^*).$$ However, pairs given by diagrams are compatible, so in particular $M_wM_w^*, M_zM_z^*$  commute (see \cite{Comp}). Hence \begin{align*}\ker M_w^*M_z^*=(\ker M_w^*\cap\ker M_z^*)\oplus \left(\ker M_w^*\cap M_w(\ker M_z^*)\right)\oplus &\left(M_z(\ker M_w^*)\cap\ker M_z^*\right)\\&\oplus \left(M_z(\ker M_w^*)\cap M_w(\ker M_z^*)\right)\end{align*}
and in turn $$C(M_w, M_z)=P_{\ker M_w^*\cap\ker M_z^*}-P_{M_w(\ker M_z^*)\cap M_z(\ker M_w^*)}.$$
\begin{remark}\label{borderprop}
    For the pair of isometries defined by the diagram $J$ we have
\begin{itemize}
\item $\ker M_w^*M_z^*=\overline{\text{Span}\{w^iz^j: (i,j)\in J_{border}\}}$
    \item $\ker M_w^*=\overline{\text{Span}\{w^iz^j: (i,j)\in J_{VB}\}}$,
    \item $\ker M_z^*=\overline{\text{Span}\{w^iz^j: (i,j)\in J_{HB}\}}$,
    \item $\ker M_w^*\cap\ker M_z^*=\overline{\text{Span}\left\{w^iz^j: (i,j)\in J_{outer}\right\}}$,
    \item $M_w(\ker M_z^*)\cap M_z(\ker M_w^*)=\overline{\text{Span}\left\{w^iz^j: (i,j)\in J_{inner}\right\}}$,
    \end{itemize}
\end{remark}
\begin{corollary}\label{defect_cor}
    The defect operator of a pair of isometries defined by the diagram $J$ is non-negative or non-positive if and only if $J$ is simple or equivalent to $\mathbb{Z}^2\setminus\mathbb{Z}_-^2$. For all other diagrams the defect is the difference of two nontrivial, mutually orthogonal projections.
\end{corollary}
\begin{proof}
    Since $C(M_w, M_z)=P_{\ker M_w^*\cap\ker M_z^*}-P_{M_w(\ker M_z^*)\cap M_z(\ker M_w^*)}$ by Remark \ref{borderprop} it is enough to show that the only diagrams for which at least one of $J_{outer}, J_{inner}$ vanish are simple diagrams or diagrams equivalent to $\mathbb{Z}^2\setminus\mathbb{Z}_-^2$.  Assume a diagram for which at least one of $J_{inner}, J_{outer}$ vanish. If there are two or more points in $J_{outer},$ then between them there is a point in $J_{inner}$ and similarly otherwise, so it is not the case. If there is precisely one point in $J_{outer}$ and $J_{inner}=\emptyset,$ then the diagram is equivalent to $\mathbb{Z}_+^2,$ so it is simple.  If there is precisely one point in $J_{inner}$ and $J_{outer}=\emptyset,$ then the diagram is equivalent to  $\mathbb{Z}^2\setminus\mathbb{Z}_-^2.$  If $J_{inner}=J_{outer}=\emptyset,$ then the diagram is equivalent to one of simple diagrams $\mathbb{Z}^2, \mathbb{Z}\times\mathbb{Z}_+,\mathbb{Z}_+\times\mathbb{Z}$. Hence $J_{inner}=\emptyset$ if and only if the diagram is simple and the only non-simple diagram for which $J_{outer}=\emptyset$ is a translation of $\mathbb{Z}^2\setminus\mathbb{Z}_-^2$. 
    \end{proof}

In the further part we will show that the shape of a diagram $J$ determines $\sigma_T(M_w, M_z)$ for $M_w, M_z\in B(H_J).$ More precisely, it is described by the parameters $\delta_\pm,\eta_\pm,\rho_\pm$ relaying on $J$ defined as follows:
\begin{itemize}
\item If $M_j=\infty$ for some $j,$ then  $\rho_-=\delta_-=\eta_-=\infty$.  If there is no $M_j=\infty,$ then (excluding the diagram $\mathbb{Z}^2$ where all $M_j=-\infty$) $M_j\in\mathbb{Z}$ for all $j\leq j_0$ for some $j_0.$ Since translation of diagrams define the same pair of isometries we may assume $j_0=0$ and define $\delta_-, \eta_-,\rho_-$ as in \eqref{deltarho+-} below. 
\item If $M_j=-\infty$ for some $j,$ then $\rho_+=\delta_+=\eta_+=\infty$. If there is no $M_j=-\infty,$ then $M_j\in\mathbb{Z}$ for all $j>j_1$ for some $j_1.$ As above we may replace $j_1$ by $-1$ and define $\delta_+, \eta_+,\rho_+$ as in \eqref{deltarho+-}.
\end{itemize}

\begin{align}\nonumber
\rho_-:=\lim_{n\to\infty}\sup_{j\leq 0} \frac{M_{j-n}-M_j}{n},\quad &\rho_+:=\lim_{n\to\infty}\sup_{j\geq -1} \frac{M_{j}-M_{j+n}}{n},\\\label{deltarho+-}\delta_-:=\lim_{n\to\infty}\inf_{j\leq 0} \frac{M_{j-n}-M_j}{n}, 
    \quad &\delta_+:=\lim_{n\to\infty}\inf_{j\geq -1} \frac{M_{j}-M_{j+n}}{n},\\\nonumber\eta_-:=\limsup_{n\to -\infty} \frac{M_0-M_n}{n}, \quad  &\eta_+:=\limsup_{n\to \infty} \frac{M_0-M_n}{n}.\end{align}

 Note that $\delta_\pm,\eta_\pm,\rho_\pm\in [0,\infty]$ and $\infty$ is possible also in \eqref{deltarho+-}.

\section{Taylor spectrum of pairs of isometries - examples}
The cases considered in \cite{Bara} yields Taylor spectrum which is not proper subset of $\cD^2$ with respect to the Lebesgue measure. More precisely, the Taylor spectrum is either the whole $\cD^2$ or the set of Lebesgue measure $0$.  Indeed, the case of non-negative defect is a doubly commuting case. Hence, by the Wold-Słociński decomposition either at least one of operators is unitary, which yields Taylor spectrum of Lebesgue measure zero or we get the pair of doubly commuting unilateral shift (which is a pair given by the diagram $\mathbb{Z}_+^2$, i.e. $H_{\ZZ_+^2}=H^2(\mathbb{T})$ - Hardy space) which has Taylor spectrum $\cD^2$. In the case of non-positive defect we may additionally get the modified bi-shift (which is a pair given by the diagram $\mathbb{Z}^2\setminus\mathbb{Z}_-^2$) which has Taylor spectrum $\cD^2$. In \cite[sec 6.2.1]{Bara} there is showed the following result. 
\begin{remark}\label{ExUnotId}
Consider the pair $(M_w^n\otimes U^k, M_w^m\otimes U^l)\in B(H^2(\mathbb{T}^2)\otimes \mathcal{H}),$ where $U\in B(\mathcal{H})$ is unitary. By \cite{Bara} we get $\sigma_T(M_w\otimes I,I\otimes U)=\cD\times \sigma(U).$ Let $$f:\mathbb{C}^2\ni (w,z)\longrightarrow (z^kw^n,z^lw^m)\in\mathbb{C}.$$

By the Spectral Mapping Theorem
$$\sigma_T(M_w^n\otimes U^k, M_w^m\otimes U^l)=\sigma(f(M_w\otimes I,I\otimes U))=f(\sigma_T(M_w\otimes I,I\otimes U))=\{(z^kw^n,z^lw^m):w\in\cD, z\in\sigma(U)\}.$$

In particular, if $U$ is a bilateral shift, we get
$$\sigma_T(M_w^n\otimes U^k, M_w^m\otimes U^l)=\{(w,z)\in\cD^2:|w|^n=|w|^m\}.$$
\end{remark}
\begin{remark}\label{ex_gp}
Let us now consider the diagram $J=\{(i,j)\in\mathbb{Z}^2:j\ge-\frac{n}{m}i\}$ and $M_w, M_z\in B(H_J)$. We may assume $m, n$ to be relatively prime and by \cite[Remark 5.1]{bur} get unique integers $0<k<n, 0\leq l<m$ such that $mk-ln=1$. Note that $J+\alpha(m,-n)\in J$ for any $\alpha\in\mathbb{Z}$ and $J+\beta(-l,k)\in J$ for $\beta\in\mathbb{N}$. Hence $U=M_{z}^{*n}M_w^m$ is unitary (precisely a bilateral shift) and $V=M_w^{*l}M_z^k$ is an isometry. Thus, by \cite[Remark 4.1]{bur}, $U$ and $V$ commute and $M_w=U^kV^n, M_z=U^lV^m.$ Since $M_w, M_z$ are unilateral shifts $V$ is a unilateral shift by \cite[Theorem 4.5]{bur}. Consequently, we may identify $H_J$ with $H^2(\mathbb{T})\otimes \ker V^*$ and $(V,U)$ with $(M_z\otimes I, I\otimes \mathcal{U})$ for $\mathcal{U}=U|_{\ker V^*}\in B(\ker V^*)$ and in turn $(M_w,M_z)$ with $(M_z^n\otimes \mathcal{U}^k, M_z^m\otimes \mathcal{U}^l).$ Since $\mathcal{U}$ is a bilateral shift, by the previous remark we get
$$\sigma_T(M_w,M_z)=\{(w,z)\in\cD^2:|w|^n=|w|^m\}.$$
\end{remark}

We may sum up the examples above as follows: if the diagram is defined by a line in the sense of Remark \ref{ex_gp}, we get Taylor spectrum of Lebesgue measure $0$ while for border defined by two rectangular half-lines (that is $\mathbb{Z}_+^2$ or $\mathbb{Z}^2\setminus\mathbb{Z}_-^2$) the spectrum is the whole $\cD.$ The main result of the paper describes the case between this two. 

The result for doubly commuting unilateral shifts extends for pairs given by diagrams contained in $\mathbb{Z}_+^2$. 
\begin{theorem}
    Let $J$ be a diagram equivalent to a subset of $\ZZ_+^2$ that is $(i,j)+J\subset \ZZ_+^2$ for some $(i,j)\in\ZZ^2$. 
    Then $$\sigma_T(M_w,M_z)=\cD^2,$$
    where $M_w,M_z\in B(H_J).$ More precisely, $\mathbb{D}^2\subset\Gamma_3(M_w, M_z).$
\end{theorem}
\begin{proof}
    Since $(i,j)+J\subset \ZZ_+^2$ we get $\sum_{(i,j)\in J}|\overline{\lambda}|^{2i}|\overline{\mu}|^{2j}< \infty$ for any $0<|\mu|<1, 0<|\lambda|<1.$ Hence $\sum_{(i,j)\in J}\overline{\lambda}^i\overline{\mu}^jw^iz^j$ is a well defined, non-trivial vector. One can check that such vector belongs to the subspace $\ker(\overline{\mu}-M_w^*)\cap \ker(\overline{\lambda}-M_z^*)$ which thus is non-empty and in turn \ref{T3} is not satisfied.

    Summing up, we get $(\mathbb{D}^2\setminus\{0\})^2\subset\Gamma_3(M_w, M_z)\subset \sigma_T(M_w,M_z)\subset\cD^2$ which by closeness of the spectrum yields the main part.

    For $\mu=0$ and $\lambda\neq 0$ we also get $\ker(M_w^*)\cap \ker(\overline{\lambda}-M_z^*)\ne\{0\}$ but vectors therein have other form. Indeed, since $(i,j)+J\subset\mathbb{Z}_+^2$ there is $(i_0,j_0)$ such that $(i_0,j)\in J_{VB}$ for any $j\ge j_0$ and such $j_0$ may be chosen minimal. Then $\sum_{j=j_0}^\infty\overline{\lambda}^jw^{i_0}z^j\in \ker M_w^*\cap \ker(\overline{\lambda}-M_z^*).$ The case $\mu\neq 0, \lambda=0$ is similar, while for $\mu=\lambda=0$ we get $\ker M_w^*\cap\ker M_z^*\neq\{0\}$ by Remark \ref{borderprop} as $J_{outer}\neq\emptyset$ for the considered diagrams.
\end{proof} 

Let us finish the section with approximations of the Taylor spectrum in the considered case.

By \cite{Wrobel} the convex hull of the Taylor spectrum is contained in the closure of the joint numerical range, i.e. 
$$conv\left(\sigma_T(T_1,T_2)\right)\subset \overline{W(T_1,T_2)},$$
where $W(T_1,T_2)=\{(\langle T_1x,x\rangle,\langle T_2x,x\rangle):x\in\mathcal{H} \textnormal{ such that }\|x\|=1\}$.

Let us consider the diagram $J.$ Since translation equivalent diagrams define the same pair of isometries we may assume without loss of generality that $\NN^2\subset J$. Let $h=\sum_{i,j\in\NN}\mu^i\lambda^jw^iz^j\in H_J$ for $\mu,\lambda\in\mathbb{D}$ and $h_1=\frac{h}{\|h\|}$ note that $\|h\|^2=\frac{1}{(1-|\mu|^2)(1-|\lambda|^2)}$. Then we have 
\begin{align*}\left(\left< M_wh_1,h_1\right>,\left< M_zh_1,h_1\right>\right)&=\frac{\left(\langle M_w h,h\rangle,\langle M_z h,h\rangle\right)}{\|h\|^2}=\\=(1-|\mu|^2)(1-|\lambda|^2)&\left(\left< \sum_{i,j\in\NN}\mu^i\lambda^jw^{i+1}z^j,\sum_{i,j\in\NN}\mu^i\lambda^jw^iz^j\right>,\left< \sum_{i,j\in\NN}\mu^i\lambda^jw^iz^{j+1},\sum_{i,j\in\NN}\mu^i\lambda^jw^iz^j\right>\right)\\=(\bar\mu,\bar\lambda)&\end{align*}

Hence $(\bar\mu,\bar\lambda)\in W(M_w,M_z)$ and $\overline{W(M_w,M_z)}=\overline{\mathbb{D}}^2$.

Similarly, by Corollarly 4.1 from \cite{KovalP} the Taylor spectrum is contained in the set $$\{(\mu,\lambda)\in \CC^2: \textnormal{ the operator } (T_1-\mu)+(T_2-\lambda)z_0 \textnormal{ is non-invertible for any } z_0\in \CC \}.$$

It requires much more calculation, but it can be shown that this result   restricts Taylor spectrum of the pair given by a diagram to the whole $\overline{\mathbb{D}}^2$ as well.  To give some hint of the proof, since we may use translation equivalent diagrams, we may assume that $(0,0)\in J$, but $(0,-1)\not\in J.$  Hence, via some laborious calculation, one can show that  for any $(\mu,\lambda)\in\mathbb{D}^2$ and each $z_0\in\CC$
there is no $h\in H_J$ such that $(M_w-\mu)h+(M_z-\lambda)z_0h=1$.

\section{Parts of Taylor spectrum of the pair of isometries given by diagrams}
In this section we describe the sets $\Gamma_i$ for the pair given by a diagram. Since the unitary part of an isometry in such pairs is a bilateral shift or vanish, its point spectrum is empty. Hence, by Remark \ref{T2iso}, $\Gamma_1(M_w, M_z)=\emptyset$ and \ref{T2} can be replaced by \eqref{T2'}.
Moreover, simple diagrams defines doubly commuting pairs which Taylor spectrum is known (see f.e. \cite{Bara}). Hence we focus on non-simple diagrams.  Since $\sigma_T(M_w,M_z)$ is closed it is enough to describe its parts up to their closure. Therefore, Theorems \ref{lem2cond} and \ref{lem3cond}, where we describe the sets $\Gamma_2(M_w,M_z)$ and $\Gamma_3(M_w,M_z)$ up to the set $\mathbb{T}^2\cup\{(0,0)\},$ are sufficient to get $\sigma_T(M_w,M_z)$. Indeed, we will see that $\mathbb{T}^2$ is in the closure of the obtained sets while $(0,0)$ is considered below. 

\begin{proposition}\label{00diagram}
    The origin $(0,0)$ belongs to $\sigma_T(M_w, M_z)$ for any $M_w, M_z\in B(H_J)$ where the diagram $J$ is non-simple. More precisely, $(0,0)\in\Gamma_2(M_w, M_z)\cap\Gamma_3(M_w, M_z)$ if the diagram is not equivalent to $\mathbb{Z}^2\setminus\mathbb{Z}_-^2$ (and non-simple), while $(0,0)\in\Gamma_2(M_w, M_z)\setminus\Gamma_3(M_w, M_z)$ for the diagram equivalent to $\mathbb{Z}^2\setminus\mathbb{Z}_-^2$.
\end{proposition}
\begin{proof}
By Corollary \ref{defect_cor} the non-simple diagram is either equivalent to $\mathbb{Z}^2\setminus\mathbb{Z}_-^2$ or its defect is a difference of two mutually orthogonal projections. 
By Remark \ref{borderprop}, similarly like in the proof of Corollary \ref{defect_cor}, we get $J_{inner}\neq\emptyset$ for all non-simple diagrams and $J_{outer}\neq\emptyset$ for all non-simple diagrams except of diagrams equivalent to $\mathbb{Z}^2\setminus\mathbb{Z}_-^2$.

By Remark \ref{borderprop} $J_{outer}\neq\emptyset$ yields \ref{T3} is not satisfied. In the case $\mathbb{Z}^2\setminus\mathbb{Z}_-^2$ note that 
$H_{\mathbb{Z}^2\setminus\mathbb{Z}_-^2}=H_{\mathbb{N}\times\mathbb{Z}_-}\oplus H_{\mathbb{Z}\times\mathbb{N}}.$  Since $H_{\mathbb{Z}\times\mathbb{N}}\subset\ran(M_w)$ as $H_{\mathbb{Z}\times\mathbb{N}}$ is the unitary subspace in the Wold decomposition of $M_w$ and $H_{\mathbb{N}\times\mathbb{Z}_-}\subset H_{\mathbb{N}\times\mathbb{Z}}\subset\ran(M_z)$ as $H_{\mathbb{N}\times\mathbb{Z}}$ is the unitary subspace in the Wold decomposition of $M_z$ we get $H_{\mathbb{Z}^2\setminus\mathbb{Z}_-^2}\subset\ran(M_w)+\ran(M_z),$ so \ref{T3} is satisfied.

To show that \ref{T2} is not satisfied take $w^iz^j$ for $(i,j)\in J_{inner}.$ By Remark \ref{borderprop} we get $w^iz^j\in M_z(\ker M_w^*)\cap M_w(\ker M_z^*)$ and in turn $h_1=M_z^* w^iz^j=w^iz^{j-1}\in\ker M_w^*, h_2=M_w^* w^iz^j=w^{i-1}z^j\in\ker M_z^*.$ On the other hand $M_zh_1-M_wh_2=0,$ so \ref{T2} is not satisfied.
\end{proof}

\subsection{Condition \ref{T2}}
The following result describes $\Gamma_2(M_w,M_z)$ with respect to the set of Lebesgue measure zero. We already explained, that such result is sufficient to describe $\sigma_T(M_w,M_z)$.

\begin{theorem}\label{lem2cond}
Let $M_w, M_z\in B(H_J)$, where $J$ is non-simple.
Then
\begin{align*}\{(\mu,\lambda)\in (\mathbb{D}\setminus\{0\})^2 : &|\mu|^{\eta_+}< |\lambda|< |\mu|^{\eta_-} \}\cup\{(0,0)\}\subset\Gamma_2(M_w,M_z)\setminus \mathbb{T}^2 \\&\subset\{(\mu,\lambda)\in (\mathbb{D}\setminus\{0\})^2 : |\mu|^{\eta_+}\leq |\lambda|\leq |\mu|^{\eta_-} \}\cup(\{0\}\times\mathbb{D})\cup(\mathbb{D}\times\{0\}),\end{align*}
where $\eta_-$ and $\eta_+$ are defined by \eqref{deltarho+-}.

Moreover,
\begin{itemize}
    \item $\mathbb{D}\times\{0\}\subset \Gamma_2(M_w,M_z)$ if and only if $M_w$ has non-trivial Wold decomposition, otherwise $\mathbb{D}\times\{0\}\cap \Gamma_2(M_w,M_z)=\{(0,0)\},$

\item $\{0\}\times\mathbb{D}\subset \Gamma_2(M_w,M_z)$ if and only if $M_z$ has non-trivial Wold decomposition,
otherwise $\{0\}\times\mathbb{D}\cap \Gamma_2(M_w,M_z)=\{(0,0)\}.$
\end{itemize}

\end{theorem}

\begin{proof}
Let us first summarize the proof, then we show its distinguished parts. 
The proof is divided among several cases, where the trivial equivalence \begin{equation}\label{equiv}(\mu,\lambda)\in\Gamma_2(M_w,M_z) \iff (\lambda,\mu)\in\Gamma_2(M_z,M_w)\end{equation} allows us to reduce their number.  
By Proposition \ref{00diagram} $(0,0)\in \Gamma_2(M_w, M_z).$
The Case II.1 shows the main part, that is  \begin{align*}\{(\mu,\lambda)\in (\mathbb{D}\setminus\{0\})^2 : |\mu|^{\eta_+}< |\lambda|< |\mu|^{\eta_-} \}&\subset\Gamma_2(M_w,M_z)\cap (\mathbb{D}\setminus\{0\})^2\\&\subset\{(\mu,\lambda)\in (\mathbb{D}\setminus\{0\})^2 : |\mu|^{\eta_+}\leq |\lambda|\leq |\mu|^{\eta_-} \}.\end{align*}
Hence we get the first inclusion. 

The second inclusion follows from the above and Case I and Case II.2 where by the last two $\mathbb{D}\times\mathbb{T}\cap\Gamma_2(M_w, M_z)=\emptyset,$ and by \eqref{equiv}   $\mathbb{T}\times\mathbb{D}\cap\Gamma_2(M_w, M_z)=\emptyset.$

The moreover part follows from  Case II.3 and Proposition \ref{00diagram} and \eqref{equiv}. 

Let us now make some preparatory work. The condition \eqref{T2'} is not satisfied if there is $h_1\notin \ran(\mu-M_w)$ such that $(\lambda-M_z)h_1\in\ran{(\mu-M_w)}.$ In all the cases below $|\mu|<1.$ Hence $\ran{(\mu-M_w)}$ is closed by Remark \ref{RS_is_cl} and so $h_1=(\mu-M_w)f+h,$ where $0\neq h\in\ker(\overline{\mu}-M_w^*).$
Clearly $(\lambda-M_z)h=(\lambda-M_z)h_1-(\lambda-M_z)(\mu-M_w)f\in\ran{(\mu-M_w)}.$
In other words \eqref{T2'} is not satisfied if and only if there is $$0\neq h\in \ker(\overline{\mu}-M_w^*),$$ such that $$(\lambda-M_z)h\in\ran (\mu-M_w).$$ 

The idea of the proof is to make attempt to construct $h$ as above in each of the considered cases. If such $h$ can be constructed, then the corresponding pair $(\mu,\lambda)\in\Gamma_2(M_w, M_z).$ Otherwise, that is if we show that such $h$ does not exist, then $(\mu,\lambda)\notin\Gamma_2(M_w, M_z).$

\textbf{Case I}: In the case we assume $\mu= 0,\lambda\in\mathbb{T}$ and show that $(\mu,\lambda) \not\in \Gamma_2(M_w,M_z)$, i.e. $\{0\}\times\mathbb{T}\cap\Gamma_2(M_w, M_z)=\emptyset.$ 

Since $\mu= 0$ we are looking for $h\in\ker M_w^*,$ which, by Remark \ref{borderprop} and \eqref{d1} has the form $h=\sum_{j:M_j\in\mathbb{Z}}a_jw^{M_j}z^j.$ Let us calculate \begin{align*}(\lambda-M_z)h&=\sum_{j:M_j\in\mathbb{Z}}\lambda a_jw^{M_j}z^j-\sum_{j:M_j\in\mathbb{Z}}a_jw^{M_j}z^{j+1}\\
&=\sum_{\substack {j:M_j\in\mathbb{Z}\\M_j<M_{j-1}}}\lambda a_jw^{M_j}z^j+\sum_{\substack {j:M_j\in\mathbb{Z}\\M_j=M_{j-1}}}\lambda a_jw^{M_j}z^j-\sum_{\substack{j:M_j\in\mathbb{Z}\\M_{j+1}=M_j}}a_jw^{M_j}z^{j+1}-\sum_{\substack{j:M_j\in\mathbb{Z}\\M_{j+1}<M_j}}a_jw^{M_j}z^{j+1}.
\end{align*}
Note that $$\sum_{\substack{j:M_j\in\mathbb{Z}\\M_{j+1}=M_j}}a_jw^{M_j}z^{j+1}=\sum_{\substack{j:M_j\in\mathbb{Z}\\M_{j+1}=M_j}}a_jw^{M_{j+1}}z^{j+1}=\sum_{\substack{j:M_j\in\mathbb{Z}\\M_{j}=M_{j-1}}}a_{j-1}w^{M_{j}}z^{j}$$ which yields
$$\sum_{\substack {j:M_j\in\mathbb{Z}\\M_j=M_{j-1}}}\lambda a_jw^{M_j}z^j-\sum_{\substack{j:M_j\in\mathbb{Z}\\M_{j+1}=M_j}}a_jw^{M_j}z^{j+1}=\sum_{\substack {j:M_j\in\mathbb{Z}\\M_j=M_{j-1}}}(\lambda a_j-a_{j-1})w^{M_j}z^j.$$
Summing up, we get $$(\lambda-M_z)h=\sum_{\substack {j:M_j\in\mathbb{Z}\\M_j<M_{j-1}}}\lambda a_jw^{M_j}z^j+\sum_{\substack {j:M_j\in\mathbb{Z}\\M_j=M_{j-1}}}(\lambda a_j-a_{j-1})w^{M_j}z^j-\sum_{\substack{j:M_j\in\mathbb{Z}\\M_{j+1}<M_j}}a_jw^{M_j}z^{j+1}.$$
Recall that we are looking for $h$ such that the $(\lambda-M_z)h\in\ran(M_w).$
Since the last sum is over $M_j$ satisfying $M_{j+1}<M_j$ which yields $(M_j-1,j+1)\in J,$ so $w^{M_j}z^{j+1}=M_w w^{M_j-1}z^{j+1}$ we get $\sum_{\substack{j:M_j\in\mathbb{Z}\\M_{j+1}<M_j}}a_jw^{M_j}z^{j+1}\in\ran(M_w)$.  On the other hand, $w^{M_j}z^j\in\ker M_w^*$ for any $j$ such that $M_j\in\mathbb{Z}.$ Consequently  $(\lambda-M_z)h\in\ran(M_w)$ if and only if the first two sums vanish, that is $a_j=0$ for $j$ such that  $M_j<M_{j-1}$ and $a_{j-1}=\lambda a_j$ for $j$ such that $M_j=M_{j-1}$. Hence $a_j\neq 0$ yields $M_j=M_{j-1}.$ Since $h\neq 0$ there is $j_0$ such that  $a_{j_0}\neq 0,$ so $M_{j_0}=M_{j_0-1}.$  If there is $k$ such that $M_{j_0}= M_{j_0-1}=\dots=M_{j_0-k}<M_{j_0-k-1},$ then $a_{j_0-k}=0$ and by  $a_{j-1}=\lambda a_j$ for all $j_0\leq j\leq j_0-k$ we get  $0=a_{j_0-k}=\lambda^k a_{j_0}\neq 0,$ a contradiction. Consequently $M_j=M_{j_0}$ for all $j<j_0$. On the other hand, if all $M_j$ were equal, then the diagram would be equivalent to $\mathbb{N}\times\mathbb{Z},$ which is not possible as it is non-simple. Hence we may choose $j_0$ the maximal, such that $M_j=M_{j_0}$ for all $j\leq j_0$. However, $a_j=\lambda^{j_0-j} a_{j_0}$ for all $j\leq j_0$ and $|\lambda|=1$ yields $\|h\|^2=\left\|\sum_{j=-\infty}^{j_0} \lambda^{j_0-j}a_{j_0}w^{M_j}z^j\right\|^2=\sum_{j\in\mathbb{N}}|\lambda|^{2j}|a_{j_0}|^2=\infty.$ In other words, $h$ does not exist, so $\{0\}\times\mathbb{T}\cap\Gamma_2(M_w, M_z)=\emptyset.$

\textbf{Case II}: $0<|\mu|<1, \lambda\in\cD$. 

Since $h\in\ker(\overline{\mu}-M_w^*)$ we get $h=\sum_{(i,j)\in J} a_j\overline{\mu}^iw^iz^j,$ where $\sum_{(i,j)\in J}|\mu|^{2i}|a_j|^2<\infty$. Let us fix $j\in\mathbb{Z}.$ Then $h_j:=\sum_{i:(i,j)\in J} a_{j}\overline{\mu}^iw^iz^{j}\in \ker(\overline{\mu}-M_w^*).$  If $M_{j}=-\infty,$ then $h_j=0.$ Indeed,  in this case $\|h_j\|^2=\sum_{i:(i,j)\in J}|a_{j}|^2|\mu|^{2i}<\infty$ if and only if $a_{j}= 0.$ On the other hand, since $0\neq h$ we get $a_j\neq 0$ for some $j,$ so there is at least one $M_j\in\mathbb{Z}.$ It follows also from the fact that the diagram is non-simple. Thus $$h=\sum_{j=j_0}^{j_1}h_j=\sum_{j=j_0}^{j_1}\sum_{i=M_j}^\infty a_j\overline{\mu}^iw^iz^j,$$ where $$j_0:=\inf\{j\in\ZZ: M_j\in\ZZ\}\quad \textnormal{and }j_1:=\sup\{j\in\ZZ: M_j\in\ZZ\}.$$

Note that
\begin{equation}\label{hnorm}
\|h_j\|^2=\left\|\sum_{i=M_j}^\infty a_j\overline{\mu}^iw^iz^j\right\|^2=|a_j|^2\sum_{i=M_j}^\infty|\mu|^{2i}=\frac{|a_j|^2|\mu|^{2M_j}}{1-|\mu|^2}.
\end{equation}

The second property of $h$ is that $(\lambda-M_z)h\in\ran(\mu-M_w).$ Let us show that it yields   \begin{equation}\label{eq_mainn}
    |\mu|^{2(M_{j-1}-M_j)}a_{j-1}=\lambda a_j
\end{equation}  for $j_0+1\leq j\leq j_1.$  Since $a_j=0$ for $j>j_1$ \eqref{eq_mainn} is clearly satisfied also for $j>j_1+1.$ However, it is not satisfied for $j=j_{1}+1$. If it were satisfied, it would yield $a_{j_1}=0$ and in turn all $a_j=0,$ so  $h=0$. Hence we prove \eqref{eq_mainn} and for the sake of completeness check coefficients for all $(i,j)\in J.$ In particular we show, that $(\lambda-M_z)h\in\ran(\mu-M_w)$ does not yield $a_{j_1}=0$.  

For this purpose we compare coefficients of
$$(\lambda-M_z)h=(\lambda-M_z)\left(\sum_{(i,j)\in J} a_j\overline{\mu}^iw^iz^j\right)=\sum_{(i,j)\in J_{HB}} 
    \lambda a_j\overline{\mu}^iw^iz^j+\sum_{(i,j)\in J\setminus J_{HB}} (\lambda a_j-a_{j-1})\overline{\mu}^iw^iz^j$$
with corresponding coefficients of an arbitrary vector in $\ran(\mu-M_w),$ that is
$$(\mu-M_w)\left(\sum_{(i,j)\in J} \alpha_{i,j}w^iz^j\right)=\sum_{(i,j)\in J_{VB}} \mu\alpha_{i,j}w^iz^j+\sum_{(i,j)\in J\setminus J_{VB}} (\mu\alpha_{i,j}-\alpha_{i-1,j})w^iz^j,$$
where obviously coefficient must satisfy $\sum_{(i,j)\in J}|\alpha_{i,j}|^2<\infty.$

If $j_1\in\mathbb{Z},$ then $j>j_1$ appears in the sums above. Note that $M_j=-\infty$ for all $j>j_1.$ Since $(i,j)\in J\setminus (J_{HB}\cup J_{VB})$ and $a_j=a_{j-1}=0$ for $j>j_1+1$ comparing coefficient  we get $\mu\alpha_{i,j}=\alpha_{i-1,j}.$  Hence $\sum_{i\in \mathbb{Z}}|\alpha_{i,j}|^2=\sum_{i\in\mathbb{Z}}|\mu|^{2i}|\alpha_{0,j}|^2$ which converge only if $\alpha_{0,j}=0$. Consequently $\alpha_{i,j}=0$ for $j>j_1+1.$

Let us now consider $j=j_1+1.$ If $i<M_{j_1},$ then $(i,j_1+1)\in J_{HB}\setminus J_{VB}$ and we get $\mu\alpha_{i,j_1+1}=\alpha_{i-1,j_1+1}.$ Hence $\alpha_{M_{j_1}-k-1,j_1+1}=\mu^k\alpha_{M_{j_1}-1,j_1+1}$ for $k\in\mathbb{N}$ and $\sum_{k\in\mathbb{N}}|\alpha_{M_{j_1}-k-1,j_1+1}|^2=\frac{|\alpha_{M_{j_1}-1,j_1+1}|^2}{1-|\mu|^2}<\infty.$ If $i=M_{j_1}+k$ for $k\in\mathbb{N},$ then $(i,j_1+1)\in J\setminus (J_{HB}\cup J_{VB})$ and we get $\mu\alpha_{M_{j_1}+k,j_1+1}-\alpha_{M_{j_1-1}+k-1,j_1+1}=-a_{j_1}\overline{\mu}^{M_{j_1}+k}$ by which
$$\mu^{k+1}\alpha_{M_{j_1}+k,j_1+1}=\mu^k\alpha_{M_{j_1-1}+k-1,j_1+1}-|\mu|^{2k}a_{j_1}\overline{\mu}^{M_{j_1}}.$$
By the inductive proof based on the equality above we get $$\mu^{k+1}\alpha_{M_{j_1}+k,j_1+1}=\alpha_{M_{j_1-1}-1,j_1+1}-(1+|\mu|^2+\dots+|\mu|^{2k})a_{j_1}\overline{\mu}^{M_{j_1}}=\alpha_{M_{j_1-1}-1,j_1+1}-\frac{1-|\mu|^{2k+2}}{1-|\mu|^2}a_{j_1}\overline{\mu}^{M_{j_1}}.$$
Since $\mu^{k+1}\alpha_{M_{j_1}+k,j_1+1}$ converge to $0$ as $k\to\infty$ taking such limit in the equality above we get $\alpha_{M_{j_1-1}-1,j_1+1}=\frac{1}{1-|\mu|^2}a_{j_1}\overline{\mu}^{M_{j_1}}$ and in turn
$$\mu^{k+1}\alpha_{M_{j_1}+k,j_1+1}=\frac{1}{1-|\mu|^2}a_{j_1}\overline{\mu}^{M_{j_1}}-\frac{1-|\mu|^{2k+2}}{1-|\mu|^2}a_{j_1}\overline{\mu}^{M_{j_1}}=\frac{|\mu|^{2k+2}}{1-|\mu|^2}a_{j_1}\overline{\mu}^{M_{j_1}}.$$ Consequently $\alpha_{M_{j_1}+k,j_1+1}=\frac{1}{1-|\mu|^2}a_{j_1}\overline{\mu}^{M_{j_1}+k+1}$ and the sequence $\sum_{k\in\mathbb{N}}\left|\frac{1}{1-|\mu|^2}a_{j_1}\overline{\mu}^{M_{j_1}+k+1}\right|^2=\frac{|a_{j_1}\overline{\mu}^{M_{j_1}+1}|^2}{(1-|\mu|^2)^3}<\infty$ for any $a_{j_1}.$  Since \eqref{eq_mainn} determines all $a_j$ on the basis of one of them, f.e. $a_{j_1}$ if $j_i\in\mathbb{Z}$ and assuming $a_{j_1}\neq 0$ the summability $\sum_{(i,j)\in J}|\mu|^{2i}|a_j|^2<\infty$ depends on $\mu, \lambda$ only.   Moreover, if $0\neq h$ is well defined, then $(\lambda-M_z)h=(\mu-M_w)\left(\sum_{(i,j)\in J} \alpha_{i,j}w^iz^j\right)$ yields by Remarks \ref{Clcond} and \ref{RS_is_cl} that  $\sum_{(i,j)\in J}|\alpha_{i,j}|^2<\infty.$ Hence all the problem of the existence of $h\in\ker(\overline{\mu}-M_w^*)$ such that $(\lambda-M_z)h\in\ran (\mu-M_w)$ reduces to the problem of summability $\sum_{(i,j)\in J}|\mu|^{2i}|a_j|^2<\infty$ which is conditioned by $\mu, \lambda.$

Let us now prove \eqref{eq_mainn}.
We start with $(i,j)\in J_{HB}.$ Since  $j_0+1\leq j\leq j_1,$ we have $M_j\in\mathbb{Z}$ and may describe $(i,j)\in J_{HB}$ as $(M_j+k,j)$ for $0\leq k <M_{j-1}-M_j$ - see \eqref{d2}. It is the case only if $M_j<M_{j-1}$. Let us show by induction that \begin{equation}\label{cos}\mu^{k+1}\alpha_{M_j+k,j}=(1+|\mu|^2+\dots+|\mu|^{2k})\lambda a_j \overline{\mu}^{M_j}=\frac{1-|\mu|^{2(k+1)}}{1-|\mu|^2}\lambda a_j \overline{\mu}^{M_j},\end{equation} for $0\leq k <M_{j-1}-M_j.$  Since $(M_j,j)\in J_{VB}$  we get base step $k=0$ directly $$\mu\alpha_{M_j,j}=\lambda a_j \overline{\mu}^{M_j}.$$ For inductive step $k>0$ we assume $$\mu^{k+1}\alpha_{M_j+k,j}=(1+|\mu|^2+\dots+|\mu|^{2k})\lambda a_j \overline{\mu}^{M_j}$$ and by $(M_j+k+1,j)\in J\setminus J_{VB}$ we get $$\mu\alpha_{M_j+k+1,j}-\alpha_{M_j+k,j}=\lambda a_j \overline{\mu}^{M_j+k+1}.$$ Hence
$$\mu^{k+2}\alpha_{M_j+k+1,j}=|\mu|^{2(k+1)}\lambda a_j \overline{\mu}^{M_j}+\mu^{k+1}\alpha_{M_j+k,j}$$ which by inductive assumption leads directly to \eqref{cos}.

In particular, for $k=M_{j-1}-M_j-1$, we get \begin{equation}\label{cosmax}\mu^{M_{j-1}-M_j}\alpha_{M_{j-1}-1,j}=\frac{1-|\mu|^{2(M_{j-1}-M_j)}}{1-|\mu|^2}\lambda a_j \overline{\mu}^{M_j}.\end{equation}

Next, we consider the case $(i,j)\in J\setminus J_{HB},$ so $\{(i,j)\in\mathbb{Z}^2:i\ge M_{j-1}\}$ - see \eqref{d2}. It is the case $M_{j-1}\in\mathbb{Z}$. For $M_j=M_{j-1}$ we get $(M_{j-1},j)=(M_j,j)\in J_{VB}$ and so $\mu\alpha_{M_j,j}=(\lambda a_j-a_{j-1})\overline{\mu}^{M_{j}}$ while for $M_j<M_{j-1}$ we get $\mu\alpha_{M_{j-1},j}-\alpha_{M_{j-1}-1,j}=(\lambda a_j-a_{j-1})\overline{\mu}^{M_{j-1}}$ and in turn by \eqref{cosmax}  \begin{align}\label{cosmax2}\mu^{M_{j-1}-M_j+1}\alpha_{M_{j-1},j}&=|\mu|^{2(M_{j-1}-M_j)}(\lambda a_j-a_{j-1})\overline{\mu}^{M_j}+\mu^{M_{j-1}-M_j}\alpha_{M_{j-1}-1,j}\\\nonumber &=|\mu|^{2(M_{j-1}-M_j)}(\lambda a_j-a_{j-1})\overline{\mu}^{M_j}+\frac{1-|\mu|^{2(M_{j-1}-M_j)}}{1-|\mu|^2}\lambda a_j \overline{\mu}^{M_j}.\end{align}
 Note that \eqref{cosmax2} for $M_j=M_{j-1}$ yields the proper formulae, so we use it for $M_j\leq M_{j-1}$.

 Let us consider remaining indices $(M_{j-1}+l,j)$ where $l\ge 1.$ Since such indices do not belong neither to $J_{VB}$ nor $J_{HB}$ we get $\mu\alpha_{M_{j-1}+l,j}-\alpha_{M_{j-1}+l-1,j}=(\lambda a_{j}-a_{j-1})\overline{\mu}^{M_{j-1}+l}$ which by induction yields 
 \begin{align}\label{ind_l} \mu^{l+1}\alpha_{M_{j-1}+l,j}&=(|\mu|^2+\dots+|\mu|^{2l})(\lambda a_j-a_{j-1})\overline{\mu}^{M_{j-1}}+\mu\alpha_{M_{j-1},j}\\\nonumber &=|\mu|^2\frac{1-|\mu|^{2l}}{1-|\mu|^2}(\lambda a_j-a_{j-1})\overline{\mu}^{M_{j-1}}+\mu\alpha_{M_{j-1},j}.\end{align}
 Since $\sum_{(i,j)\in J}|\alpha_{i,j}|^2<\infty$ the sequence $\alpha_{M_{j-1}+l,j}\to 0,$ as $l\to\infty,$ so taking such limit in \eqref{ind_l} we get 

$$0=|\mu|^2\frac{1}{1-|\mu|^2}(\lambda a_j-a_{j-1})\overline{\mu}^{M_{j-1}}+\mu\alpha_{M_{j-1},j}$$ which multiplied by $\mu^{M_{j-1}-M_j}$ and by \eqref{cosmax2} yields
$$0=|\mu|^{2(M_{j-1}-M_j+1)}\frac{1}{1-|\mu|^2}(\lambda a_j-a_{j-1})\overline{\mu}^{M_j}+|\mu|^{2(M_{j-1}-M_j)}(\lambda a_j-a_{j-1})\overline{\mu}^{M_j}+\frac{1-|\mu|^{2(M_{j-1}-M_j)}}{1-|\mu|^2}\lambda a_j \overline{\mu}^{M_j}.$$
Hence one can directly calculate  \eqref{eq_mainn}.

Note that if $j_0\in\ZZ$ (so $M_{j_0-1}=\infty$), we can take $k\to\infty$ in \eqref{cos} and get $a_{j_0}=0$, for $\lambda\not=0$. Then by \eqref{eq_mainn} all $a_j$ vanish as well, so $h=0.$ In other words,  for $\lambda\not=0$ the assumption  $h\neq 0$ yields $j_0=-\infty$, so $M_j\in \ZZ$ for any $j\leq j_1.$

Since equivalent diagrams  define the same pair of isometries we may assume without loss of generality that $M_0=0$ (so $j_0\leq 0\leq j_1$) and $a_0=1.$

In the considered case $\mu\neq 0$. If also $\lambda\neq 0,$ then by \eqref{eq_mainn} we get

$$a_{j-1}=\frac{\lambda}{|\mu|^{2(M_{j-1}-M_j)}}a_j\quad\textnormal{and}\quad a_{j+1}=\frac{|\mu|^{2(M_{j}-M_{j+1})}}{\lambda}a_j$$ and in turn (assuming $j_0=-\infty$)
\begin{equation}\label{ajot}a_j=\left\{\begin{array}{l}\frac{\lambda^{-j}}{|\mu|^{2(M_j-M_0)}}a_0, \quad \textnormal{ for } j\leq 0\\\frac{|\mu|^{2(M_0-M_j)}}{\lambda^j}a_0, \quad \textnormal{ for } 0\leq j\leq j_1\end{array}\right\}=\frac{|\mu|^{2(M_0-M_j)}}{\lambda^j}=\frac{|\mu|^{-2M_j}}{\lambda^j}.\end{equation}

{\bf Case II.1}: In this case, we assume $0<|\mu|<1,\; 0<|\lambda|< 1$ and show that $|\mu|^{\eta_+}\leq |\lambda|\leq |\mu|^{\eta_-}$ is the necessary condition and $|\mu|^{\eta_+}< |\lambda|< |\mu|^{\eta_-}$ is the sufficient condition for $(\mu,\lambda)\in\Gamma_2(M_w, M_z).$ 

First let us consider the case $j_0\in\mathbb{Z}$. Then, as we already showed $h$ does not exist, so none of $(\mu,\lambda)$ satisfying the assumption belongs to $\Gamma_2(M_w,M_z)$. Hence the condition we are going to prove should give empty set for $j_0\in\mathbb{Z}$. Indeed, $j_0\in\mathbb{Z}$ yields $M_{j_0-1}=\infty,$ so $\eta_-=\infty$ and the part of the necessary condition $|\lambda|\leq |\mu|^{\eta_-}=0$ is not satisfied by any $\lambda\neq 0$.
  
For $j_0=-\infty,$ by  \eqref{hnorm} and \eqref{ajot} we get
\begin{align*}\|h\|^2=\sum_{j=-\infty}^{j_1} |a_j|^2\frac{|\mu|^{2M_j}}{1-|\mu|^2}&=\sum_{j=-\infty}^{j_1} \Big|\frac{|\mu|^{-2M_j}}{\lambda^j}\Big|^2\frac{|\mu|^{2M_j}}{1-|\mu|^2}=\sum_{j=-\infty}^{j_1} \frac{|\mu|^{-2M_j}}{|\lambda|^{2j}}\frac{1}{1-|\mu|^2}\\&=\frac{1}{1-|\mu|^2}\left(\sum_{j=-\infty}^{0} \frac{|\mu|^{-2M_j}}{|\lambda|^{2j}}+\sum_{j=1}^{j_1} \frac{|\mu|^{-2M_j}}{|\lambda|^{2j}}\right).\end{align*}

Hence $h$ is well defined if and only if both the series above are convergent, where in the case $j_1\in\mathbb{Z}$ the second one is finite.
The Cauchy root test of convergence gives an answer. Indeed,
$$\limsup_{n\to -\infty}\sqrt[-n]{\frac{|\mu|^{-2M_n}}{|\lambda|^{2n}}}=\limsup_{n\to -\infty}\frac{|\lambda|^{2}}{|\mu|^{2\frac{-M_n}{n}}}=\frac{|\lambda|^{2}}{|\mu|^{2\eta_-}}\leq 1\textnormal{ is a necessary and }\frac{|\lambda|^{2}}{|\mu|^{2\eta_-}}<1 \textnormal{ is a sufficient}$$
condition for convergence of the first series and in the case $j_1=\infty$ 
$$\limsup_{n\to \infty}\sqrt[n]{\frac{|\mu|^{-2M_n}}{|\lambda|^{2n}}}=\limsup_{n\to \infty}\frac{|\mu|^{2\frac{-M_n}{n}}}{|\lambda|^{2}}=\frac{|\mu|^{2\eta_+}}{|\lambda|^2}\leq 1\textnormal{ is a necessary and } \frac{|\mu|^{2\eta_+}}{|\lambda|^2}<1 \textnormal{ is a sufficient}$$
condition for convergence of the second series. Note that for $j_1$ finite the series is finite, and we have $\eta_+=\infty$, so $\frac{|\mu|^{2\eta_+}}{|\lambda|^2}=0\leq 1.$ In other words, we can use $\frac{|\mu|^{2\eta_+}}{|\lambda|^2}\leq 1$ as a universal condition for $j_1\in\mathbb{Z}\cup\{\infty\}.$

Summing up, we have showed that if there is $h\in\ker(\overline{\mu}-M_w^*)$ such that $(\lambda-M_z)h\in\ran (\mu-M_w),$ then it is of the form $h=\sum_{j=j_0}^{j_1}h_j=\sum_{j=j_0}^{j_1}\sum_{i=M_j}^\infty a_j\overline{\mu}^iw^iz^j,$ where $a_j$ are in relation \eqref{ajot}.  By \eqref{ajot} such $h$ is well defined  if $|\mu|^{\eta_+}< |\lambda|< |\mu|^{\eta_-}$ (sufficient condition) and only if $|\mu|^{\eta_+}\leq |\lambda|\leq |\mu|^{\eta_-}$ (necessary condition).
Since existence of such $h$ is equivalent to $(\mu,\lambda)\in\Gamma_2(M_w, M_z)$ we get the sufficient and necessary conditions of $(\mu,\lambda)\in\Gamma_2(M_w, M_z)$.

{\bf Case II.2}: In this case, we assume $0<|\mu|<1,\; |\lambda|=1$ and show that $\mathbb{D}\setminus\{0\}\times\mathbb{T}\cap\Gamma_2(M_w, M_z)=\emptyset.$

The case $j_0\in\mathbb{Z}$ is clear. For $j_0=-\infty$ by \eqref{hnorm} and \eqref{ajot} we get
$$\|h\|^2=\sum_{j=-\infty}^{j_1} |a_j|^2\frac{|\mu|^{2M_j}}{1-|\mu|^2}
=\sum_{j=-\infty}^{j_1} |\mu|^{-2M_j}\frac{1}{1-|\mu|^2}.$$
However, $|\mu|<1$ and $M_j\geq M_0 = 0$ for $j<0$ yields $\|h\|=\infty,$ a contradiction. Hence $\mathbb{D}\setminus\{0\}\times\mathbb{T}\cap\Gamma_2(M_w, M_z)=\emptyset.$

{\bf Case II.3}: In this case, we assume $0<|\mu|<1,\;\lambda=0$ and show that: 

- if $M_w$ has nontrivial Wold decomposition, then $\left(\mathbb{D}\setminus\{0\}\right)\times \{0\}\subset\Gamma_2(M_w, M_z),$

- if $M_w$ is a unilateral shift, then $\left(\mathbb{D}\setminus\{0\}\right)\times \{0\}\cap\Gamma_2(M_w, M_z)=\emptyset.$

By \eqref{eq_mainn} and \eqref{cos} we have $a_j=0$ for any $j_0-1< j< j_1$.

If $j_1\in \ZZ,$ then $h=\sum_{i\geq M_{j_1}} \overline{\mu}^i w^iz^{j_1}\in \ker(\overline{\mu}-M_w^*)$. Note that the space $\overline{\textnormal{Span}\{w^iz^{j}: i\in \ZZ, j>j_1\}}$  is unitary space in Wold decomposition of $M_w$ (precisely $M_w$ is a bilateral shift there). Since $|\mu|<1$ it belongs to the resolvent of a bilateral shift, i.e. $\mu-M_w$ is invertible on the unitary part of $M_w,$ we get
$-M_zh=-\sum_{i\geq M_{j_1}} \overline{\mu}^i w^iz^{j_1+1}\in \ran(\mu-M_w)$. Hence we get the first part.

If $j_1=\infty,$ then $a_j=0$ for any $j_0-1< j< j_1$ is equivalent to $h=0$. In other words, if $M_w$ has a trivial Wold decomposition (it is a unilateral shift), then $(\mu,0)\not\in \Gamma_2(M_w,M_z)$. Hence we get the second part.

\end{proof}

\subsection{Condition \ref{T3}} Power partial isometries are partial isometries whose powers are partial isometries as well.
 The class is described in \cite{HW} where the decomposition of power partial isometry among unitary operator, unilateral shift, backward shift and truncated shifts is showed. 
Recall that a truncated shift of index $k$ is an operator $T_k\in B(\mathbb{C}^k)$ given by a matrix $$T_k=\left[\begin{matrix}
 0 & 0 & 0 &  \hdots & 0\\
 1 & 0 & 0 & \hdots & 0 \\
 0 & 1 & 0 & \hdots & 0 \\
 \vdots&& \ddots &\ddots &  \vdots \\
 0 & 0 &  \hdots & 1 & 0
\end{matrix}\right].$$ 

More generally, $T_k\otimes I_{\mathcal{H}}\in B(\mathbb{C}^k\otimes \mathcal{H})$ is a truncated shift of index $k$ and multiplicity $\dim\mathcal{H}.$   
\begin{remark}\label{ts}
Let $T_k\in B(\mathbb{C}^k)$ be a truncated shift of index $k$ and $\mathcal{H}$ be a Hilbert space. Then:
 \begin{itemize}
         \item $\sigma(T_k\otimes I_{\mathcal{H}})=\sigma (T_k)=\{0\},$
     \item if $\{k_i\}_{i\ge 0}$ is a strongly increasing sequence of integers, then $\sigma\left(\bigoplus_{i\ge 0}T_{k_i}\right)=\cD$. More precisely, $\ran\left(\lambda-\bigoplus_{i\ge 0}T_{k_i}\right)$ is not closed for $\lambda\in \cD\setminus \{0\}$.
     
 \end{itemize}
\end{remark}
\begin{proof}
    The spectrum of a truncated shift of index $k$ follows directly from the matrix representation. 
    
    For the second condition, by Remark \ref{Clcond} $\ran\left(\lambda-\bigoplus_{i\ge 0}T_{k_i}\right)$ is closed if and only if $\lambda-\bigoplus_{i\ge 0}T_{k_i}$ is bounded below which in turn holds if and only if $(\lambda-T_{k_i})^{-1}$ exist for all $k_i$ and are jointly bounded. However, $$(\lambda-T_{k_i})^{-1}=\left[\begin{matrix}
 \frac{1}{\lambda} & 0 & 0 &  \hdots & 0\\
 \frac{1}{\lambda^2} & \frac{1}{\lambda} & 0 & \hdots & 0 \\
 \frac{1}{\lambda^3} & \frac{1}{\lambda^2} & \frac{1}{\lambda} & \hdots & 0 \\
 \vdots&& \ddots &\ddots &  \vdots \\
 \frac{1}{\lambda^{k_i+1}} & \frac{1}{\lambda^{k_i}} &  \hdots & \frac{1}{\lambda^2} & \frac{1}{\lambda}
\end{matrix}\right]$$ for $\lambda\neq 0$.
Hence $\|(\lambda-T_{k_i})^{-1}\|\ge \|(\lambda-T_{k_i})^{-1}(1,0,\dots,0)\|=\sqrt{\frac{1}{|\lambda|^{2(k_i+1)}}+\dots+\frac{1}{|\lambda|^2}}\geq\sqrt{k_i}.$ Since the sequence $k_i$ is strongly increasing sequence of integers, it diverges to $\infty,$ so  $\|(\lambda-T_{k_i})^{-1}\|$ are not jointly bounded.
\end{proof}
Let $r(T)$ denotes the spectral radius of an operator $T$ and  $i(T)=\lim_{n\to\infty}(\inf\{\|T^n x\|:\|x\|=1\})^{\frac1n}.$ Spectrum of weighted shifts is described in \cite{Ridge}. We recall the result in the case of non-zero weights.

\begin{theorem}\label{ridge}
 Let $B_w$ on $l^2(\mathbb{Z})$ and $S_w$ on $l^2(\mathbb{Z}_+)$ be a bilateral weighed shift and right weighed shift  respectively, with weights $w_j\in\mathbb{C}\setminus\{0\}.$ Moreover, let

  $$i(B_w)^-=\lim_{n\to\infty}\inf_{j\leq 0}|w_{j-1}\dots w_{j-n}|^\frac1n,\quad i(B_w)^+=\lim_{n\to\infty}\inf_{j\geq -1}|w_{j+1}\dots w_{j+n}|^\frac1n,$$
 $$r(B_w)^-=\lim_{n\to\infty}\sup_{j\leq 0}|w_{j-1}\dots w_{j-n}|^\frac1n,\quad r(B_w)^+=\lim_{n\to\infty}\sup_{j\geq -1}|w_{j+1}\dots w_{j+n}|^\frac1n.$$ 
 Then:
 
 \begin{itemize}
     \item $i(B_w)=\min\{i(B_w)^-, i(B_w)^+\}, r(B_w)=\max\{r(B_w)^-, r(B_w)^+\},$
     \item $i(S_w)=\lim_{n\to\infty}\inf_{j}|w_{j+1}\dots w_{j+n}|^\frac1n, r(S_w)=\lim_{n\to\infty}\sup_{j}|w_{j+1}\dots w_{j+n}|^\frac1n,$
     \item $\sigma_{ap}(S_w)=\{\lambda\in\mathbb{C}: i(S_w)\leq |\lambda|\leq r(S_w)\}$
     \item $\sigma_{ap}(S_w^*)=\sigma(S_w^*)=\sigma(S_w)=\{\lambda\in\mathbb{C}: |\lambda|\leq r(S_w)\},$
     \item $\sigma_{ap}(B_w)=\{\lambda\in\mathbb{C}: i(B_w)^+\leq |\lambda|\leq r(B_w)^+\text{ or }r(B_w)^+\leq |\lambda|\leq i(B_w)^-\text{ or }i(B_w)^-\leq |\lambda|\leq r(B_w)^-\}.$
 \end{itemize}
\end{theorem}
\begin{proof}
    As we mentioned, the Theorem follows directly from \cite{Ridge}. In particular, formulae on $\sigma_{ap}(S_w)$ is \cite[Theorem 1]{Ridge}, the part $\sigma_{ap}(S_w^*)=\sigma(S_w)=\sigma(S_w^*)$ is \cite[Theorem 6]{Ridge}, while formulae on $\sigma(S_w)$ is the first Corollary in \cite{Ridge}. Let us only comment the formulae on $\sigma_{ap}(B_w)$ which is \cite[Theorem 3]{Ridge}. However, in \cite[Theorem 3]{Ridge} the spectrum is described by two cases: $$i(B_w)^-\leq |\lambda|\leq r(B_w)^-\text{ or }i(B_w)^+\leq |\lambda|\leq r(B_w)^+ \;\text{ if }\; r(B_w)^-<i(B_w)^+$$ and $$i(B_w)\leq |\lambda|\leq r(B_w) \text{ otherwise.}$$ Let us check that it is the same as in the statement. Indeed, we have $i(B_\omega)^-\leq r(B_\omega)^-$ and $i(B_\omega)^+\leq r(B_\omega)^+$. Hence, if 
    $r(B_w)^-<i(B_w)^+,$ then $i(B_w)^-<r(B_w)^+$ and the middle term $r(B_w)^+\leq |\lambda|\leq i(B_w)^-$ is not satisfied by any $\lambda.$ The two other conditions in the statement are the same as in \cite[Theorem 3]{Ridge}. Assume now  $i(B_w)^+\leq r(B_w)^-.$ If  $r(B_w)^+\leq i(B_w)^-,$ then  the three conditions in the statement sum up to $i(B_w)^+\leq |\lambda|\leq r(B_w)^-$ and $i(B_w)^+=i(B_w)$ and $r(B_w)^-=r(B_w)$. If $r(B_w)^+> i(B_w)^-,$ then $i(B_w)^+\leq |\lambda|\leq r(B_w)^+$ and $i(B_w)^-\leq |\lambda|\leq r(B_w)^-$ has non-trivial intersection and hence their sum is $i(B_w)\leq |\lambda|\leq r(B_w).$ 
\end{proof}

The convention to assume inequalities with $0^0$ or $1^\infty$ as satisfied allows us to formulate the following result in a relatively concise way. However, the convention yields, see  Remark \ref{pab}, that the set $\{(0,0)\}\cup\mathbb{T}^2$ is always included in the obtained set while it is not necessarily in $\Gamma_3(M_w,M_z).$ In other words the result describes $\Gamma_3(M_w,M_z)$ up to the set $\{(0,0)\}\cup\mathbb{T}^2.$ We have already explained that such description is sufficient to calculate $\sigma_T(M_w,M_z).$  

\begin{theorem}\label{lem3cond}
Let $M_w, M_z\in \mathcal{B}(H_J)$, where  $J$ is a non-simple diagram  and $\delta_-,\delta_+,\rho_-,\rho_+$ be given by \eqref{deltarho+-}.

Then:
\begin{itemize}
 \item if both $M_w, M_z$ have non-trivial Wold decompositions, then
$$\Gamma_3(M_w,M_z)\cup \{(0,0)\}\cup\mathbb{T}^2=\mathbb{T}\times\overline{\mathbb{D}}\cup \overline{\mathbb{D}}\times\mathbb{T} \cup \{(0,0)\},$$ 
  \item if $M_z$ has non-trivial Wold decomposition and $M_w$ is a unilateral shift, then
    $$\Gamma_3(M_w,M_z)\cup\mathbb{T}^2=\{(\mu,\lambda)\in \overline{\mathbb{D}}^2: |\mu|^{\rho_+}\leq |\lambda|\leq |\mu|^{\delta_+}\}\cup\cD\times\mathbb{T},$$
    \item if $M_w$ has non-trivial Wold decomposition and $M_z$ is a unilateral shift, then
    $$\Gamma_3(M_w,M_z)\cup\mathbb{T}^2=\{(\mu,\lambda)\in\overline{\mathbb{D}}^2: |\mu|^{\rho_-}\leq |\lambda|\leq |\mu|^{\delta_-}\}\cup\mathbb{T}\times\cD,$$
    \item if $M_w,M_z$ are unilateral shifts, then    $$\Gamma_3(M_w,M_z)\cup\mathbb{T}^2=\{(\mu,\lambda)\in\overline{\mathbb{D}}^2:
|\mu|^{\rho_-}\leq|\lambda|\leq |\mu|^{\delta_-}   \textnormal{ or }
|\mu|^{\delta_-}\leq|\lambda|\leq |\mu|^{\rho_+}  \textnormal{ or }
|\mu|^{\rho_+}\leq|\lambda|\leq |\mu|^{\delta_+}\},$$
\end{itemize}

where inequalities with $1^\infty$ or $0^0$ are assumed to be satisfied.
\end{theorem}

\begin{proof}
 Since the diagram is non-simple, there are possible only the cases as in the statement. Let us remark that $\{(0,0)\}$ is added in the first case only because of diagrams equivalent to $\mathbb{Z}^2\setminus\mathbb{Z}_-^2$ which by Proposition \ref{00diagram} are the only (non-simple) diagrams such that $(0,0)\notin\Gamma_3(M_w,M_z)$. 
 
 The plan of the proof:
\begin{itemize}
    \item In the first step we investigate $\Gamma_3(M_w,M_z)$ on  $\mathbb{D}\setminus\{0\}\times\overline{\mathbb{D}}$ - Case I and on $\{0\}\times\mathbb{T}$ - Case II. Since the proof require different approach depending on a type of the diagram, each case is divided among four subcases.
    \item In the second step we consider symmetric cases $\overline{\mathbb{D}}\times\mathbb{D}\setminus\{0\}$ - Case I' and  $\mathbb{T}\times\{0\}$ - Case II'. We immediately get the result from Cases I and II. However, the symmetry divide diagrams among different four types.
      \item In the final step we combine the types of diagrams from Cases I, II, with the types from Cases I', II'. Fortunately, not all combinations are possible, and some are symmetric to the other. Hence we need to show the result only for 6 types of diagrams: Cases A-F.
\end{itemize} 
 {\bf Step I}
 
 Since the diagram is non-simple there is at least one $M_j\in\mathbb{Z}$ and we may define
 $$j_0:=\inf\{j\in\ZZ: M_j\in\ZZ\} \quad \text{and}\quad j_1:=\sup\{j\in\ZZ: M_j\in\ZZ\},$$
    where $-\infty\leq j_0\leq j_1\leq\infty$. The diagrams are divided according to finiteness of $j_0, j_1$.
    
    Moreover, denote $H_j:=\bigoplus_{i\geq M_j} \mathbb{C}w^{i}z^{j}$, for $j_0-1<j<j_1+1$ and $H_{uw}=H_{\mathbb{Z}\times\{j:j>j_1\}}$ if $j_1\in\mathbb{Z}$ and $H_{uw}=\{0\}$ if $j_1=\infty$. Thus $H_J=H_{uw}\oplus \bigoplus_{j=j_0}^{j_1}H_j$, where $H_{uw}$ is the unitary subspace in the Wold decomposition of $M_w$.

{\bf Case I }: $0<|\mu|<1, \lambda\in\overline{\mathbb{D}}$. 

It is clear that $H_j$ reduces $M_w$ to a unilateral shift and for $e_j:=\sum_{i=M_j}^\infty \overline{\mu}^{i-M_j}w^iz^{j}$ we get $(\mu-M_w)H_j=H_j\ominus \mathbb{C}e_j$, where the closeness of the left hand side we get by Remark \ref{RS_is_cl}. Hence $$ \bigoplus_{j=j_0}^{j_1} (H_j\ominus \mathbb{C}e_j)\subset \ran(\mu-M_w).$$
Since $\bigoplus_{j=j_0}^{j_1}\mathbb{C}e_j\subset\ran(\mu-M_w)^\perp$ the decomposition $H=H_{uw}\oplus \bigoplus_{j=j_0}^{j_1} (H_j\ominus \mathbb{C}e_j)\oplus \bigoplus_{j=j_0}^{j_1}\mathbb{C}e_j,$ yields $\ran(\mu-M_w)=H_{uw}\oplus \bigoplus_{j=j_0}^{j_1} (H_j\ominus \mathbb{C}e_j)$ and   $\ran(\mu-M_w)^\perp=\bigoplus_{j=j_0}^{j_1}\mathbb{C}e_j.$  Hence it is enough to check whether $\bigoplus_{j=j_0}^{j_1}\mathbb{C}e_j\subset \ran(\mu-M_w)+\ran(\lambda-M_z).$ By the closeness of $\ran(\mu-M_w)$ we get $$\ran(\mu-M_w)+\ran(\lambda-M_z)=\ran(\mu-M_w)\oplus P_{\ran(\mu-M_w)^\perp}\ran(\lambda-M_z).$$ However, by commutativity of $M_w, M_z$ we get $P_{\ran(\mu-M_w)^\perp}\ran(\lambda-M_z)=P_{\ran(\mu-M_w)^\perp}(\lambda-M_z)\ran(\mu-M_w)^\perp.$ Recall that $F:=P_{\ran(\mu-M_w)^\perp}(\lambda-M_z)\in B(\ran(\mu-M_w)^\perp)$ is a fringe operator defined in \cite{Y}. Let us calculate $F$ precisely. Since for $j_1\in \ZZ$, $M_z H_{j_1}\subset H_{uw}\subset\ran(\mu-M_w)$ we get $P_{\ran(\mu-M_w)^\perp}(\lambda-M_z)e_{j_1}=\lambda e_{j_1}.$ 
On the other hand,  $$P_{\ran(\mu-M_w)^\perp}(\lambda-M_z)e_j=\lambda e_j - \mu^{M_j-M_{j+1}}e_{j+1},$$ for $j_0-1<j< j_1.$ Indeed, $M_z H_{j}\subset H_{j+1}$ yields $(\lambda-M_z)(e_j)=\lambda e_j+\alpha_je_{j+1}+a_{j+1}$, where $a_{j+1}\in H_{j+1}\ominus \mathbb{C}e_{j+1}\subset \ran(\mu-M_w).$ We get $\alpha_j=-\mu^{M_j-M_{j+1}}$ by  
\begin{align*}\alpha_j\|e_{j+1}\|^2&=\langle (\lambda-M_z)e_j, e_{j+1}\rangle=\langle -M_ze_j, e_{j+1}\rangle=\langle -\sum_{i\geq M_{j}} \overline{\mu}^{i-M_j} w^{i}z^{j+1}, \sum_{i\geq M_{j+1}} \overline{\mu}^{i-M_{j+1}} w^{i}z^{j+1} \rangle\\&=-\mu^{M_j-M_{j+1}}\|\sum_{i\geq M_{j}} \overline{\mu}^{i-M_j} w^{i}z^{j+1}\|^2=-\mu^{M_j-M_{j+1}}\|e_{j}\|^2=-\mu^{M_j-M_{j+1}}\|e_{j+1}\|^2,\end{align*} where 
$\|e_{j}\|^2=\|e_{j+1}\|^2$ as they both are equal to $\sum_{i=0}^\infty|\mu^i|^2=\frac{1}{1-|\mu|^2}.$

Summing up,  
$$Fe_j=\begin{cases}
    \lambda e_j-\mu^{M_j-M_{j+1}}e_{j+1}, &\textnormal{ for } j_0-1<j< j_1,\\
 \lambda e_j,   \textnormal{ for } j=j_1, &\textnormal{ if }j_1\in\mathbb{Z}
\end{cases}$$ and 

\begin{equation}\label{war}
H_{uw}\oplus\bigoplus_{j=j_0}^{j_1} (H_j\ominus \mathbb{C}e_j)\oplus \ran(F) =\ran(\mu-M_w)+\ran(\lambda-M_z).  
\end{equation}
Hence $(\mu,\lambda)\in\Gamma_3(M_w, M_z)$ if and only if $\ran(F)\neq\ran(\mu-M_w)^\perp,$ in particular if $\ran(F)$ is not closed.  

{\bf Case I.1}: $j_0=-\infty$, $j_1=\infty,$ so all $M_j$ are finite.
    
For such diagrams $F=\lambda -B_\omega$, where $B_\omega$ is a bilateral shift with the weights $\omega_j:=\mu^{M_j-M_{j+1}}$. Hence $\ran(F)\neq\ran(\mu-M_w)^\perp$ 
if and only if $\overline{\lambda}\in\sigma_{ap}(B^*_\omega)$ (see \cite[chap. XI, Proposition 1.1]{Conway}).

Note that $B_{\omega}^*$ is a bilateral shift with weights $\{\overline{\omega}_{-j-1}\}_{j\in\ZZ}$. Hence taking advantage of Theorem \ref{ridge} we may calculate \begin{align*}i(B_\omega^*)^-&=i(B_\omega)^+=\lim_{n\to\infty}\inf_{j\geq -1} |\omega_j\omega_{j+1}\cdot\ldots\cdot\omega_{n+j-1}|^\frac1n\\&=\lim_{n\to\infty}\inf_{j\geq -1} |\mu^{M_j-M_{j+1}}\mu^{M_{j+1}-M_{j+2}}\cdot\ldots\cdot\mu^{M_{n+j-1}-M_{n+j}}|^\frac1n=|\mu|^{\lim_{n\to\infty}\sup_{j\geq -1} \frac{M_j-M_{j+n}}{n}}=|\mu|^{\rho_+}.\end{align*}

In a similar way one may check that $i(B_w^*)^+=|\mu|^{\rho_-}, r(B_w^*)^-=|\mu|^{\delta_+}, r(B_w^*)^+=|\mu|^{\delta_-}.$

By Theorem \ref{ridge}
$(\mu,\lambda)\in \Gamma_3(M_w,M_z)$ if and only if $|\mu|^{\rho_-}\leq|\lambda|\leq |\mu|^{\delta_-}$ or $|\mu|^{\rho_+}\leq|\lambda|\leq |\mu|^{\delta_+}$ or $|\mu|^{\delta_-}\leq|\lambda|\leq |\mu|^{\rho_+}.$

{\bf Case I.2}: $j_0=-\infty$ and $j_1\in\ZZ$. 

Note that $j_1\in\mathbb{Z}$ yields $\rho_+=\delta_+=\infty.$

In this case $F=\lambda Id-S^*_{\omega}$, where $S_\omega$ is a unilateral shift with  weights $\omega_j:=\overline{\mu}^{M_{j-1}-M_{j}}$, for $j\leq j_1$.

Again $(\mu,\lambda)\in \Gamma_3(M_w,M_z)$ if and only if $\overline{\lambda}\in \sigma_{ap}(S_\omega)$. We calculate $i(S_\omega)=|\mu|^{\rho^-}, r(S_\omega)=|\mu|^{\delta^-},$ as in the previous case. Hence, by Theorem \ref{ridge}, 
$(\mu,\lambda)\in \Gamma_3(M_w,M_z)$ if and only if $ |\mu|^{\rho_-}\leq|\lambda|\leq |\mu|^{\delta_-}$.

{\bf Case I.3}: $j_0\in\ZZ$ and $j_1=\infty$. 

Since $j_0\in\mathbb{Z}$ we get $\delta_-=\rho_-=\infty.$ Moreover, $F=\lambda -S_{\omega}$, where $S_\omega$ is a unilateral shift with the weights $\omega_j:=\mu^{M_j-M_{j+1}}$.
Again $(\mu,\lambda)\in \Gamma_3(M_w,M_z)$ if and only if $\overline{\lambda}\in \sigma_{ap}(S^*_\omega)$. We calculate $r(S_w)=|\mu|^{\delta_+}$ as previously and by Theorem \ref{ridge} get $(\mu,\lambda)\in \Gamma_3(M_w,M_z)$ if and only if $|\lambda|\leq |\mu|^{\delta_+}$.

{\bf Case I.4}: $j_0,j_1\in\ZZ$.

In this case  $\delta_-=\rho_-=\delta_+=\rho_+=\infty.$ Moreover, $\dim(\ran(\mu-M_w)^\perp)<\infty$ and $$Fe_j=\begin{cases}
    -\mu^{M_j-M_{j+1}}e_{j+1}+\lambda e_j, & \textnormal{ for } j=j_0,j_0+1,\ldots,j_1-1,\\
    \lambda e_j, & \textnormal{ for } j=j_1.
\end{cases}$$
So, $\ran(F)=\ran(\mu-M_w)^\perp$ if and only if $\lambda\not=0$. Thus $(\mu,\lambda)\in \Gamma_3(M_w,M_z)$ if and only if $\lambda=0$.

{\bf Case II }: $\mu=0$ and $|\lambda|=1$.

 We proceed similarly like in the previous case, but we get different fringe operator. Instead of $e_j$ we take $f_j:=w^{M_j}z^{j}\in H_j\cap \ker(M_w^*)$, for $j_0-1<j<j_1+1$. Since $M_w|_{H_j}$ is a unilateral shift we have $$M_w(H_j)=H_j\ominus \mathbb{C}f_j,\textnormal{ for } j_0-1<j<j_1+1$$
    and $$M_w(H_{uw})=H_{uw}.$$     
   Note that, if $M_{j}=M_{j+1}$, then $M_z(H_j)= H_{j+1}$ while for $ M_{j+1}<M_j$ we get $M_z(H_j)\subset H_{j+1}\ominus \mathbb{C}f_{j+1}.$
    Indeed, $w^iz^{j+1}\perp f_{j+1}$ for $i> M_{j+1}.$ Hence, if $M_j>M_{j+1},$ then $M_z\left(\sum_{i\geq M_{j}} \alpha_i w^iz^{j}\right)=\sum_{i\geq M_{j}} \alpha_i w^iz^{j+1}\perp f_{j+1},$ where $\sum_{i\geq M_{j}} \alpha_i w^iz^{j}$ is an arbitrary vector in $H_j.$ In particular, $P_{\ran(M_w)^\perp}(\lambda-M_z)f_j=\lambda f_j$ if $j=j_1\in\mathbb{Z}$ or $M_j>M_{j+1}$ and $P_{\ran(M_w)^\perp}(\lambda-M_z)f_j=\lambda f_j-f_{j+1}$ if $M_j=M_{j+1}.$

Therefore, 
\begin{equation}\label{war2}
H_{uw}\oplus\bigoplus_{j=j_0}^{j_1} (H_j\ominus \mathbb{C}f_j)\oplus \ran(F) =\ran(M_w)+\ran(\lambda-M_z).  
\end{equation}
where $F\in B(\ran(M_w)^\perp)$ is a fringe operator $F=P_{\ran(M_w)^\perp}(\lambda-M_z)$ where 
$$Ff_j=\begin{cases}
    \lambda f_j,&\textnormal{ for } j=j_1\textnormal{ if }j_1\in\mathbb{Z},\\
    \lambda f_j-f_{j+1},&\textnormal{ for }j_0-1<j<j_1 \textnormal{ and } M_j=M_{j+1},\\
    \lambda f_j,& \textnormal{ for }j_0-1<j<j_1 \textnormal{ and } M_j>M_{j+1}.
\end{cases}$$
Again, we need to check whether $\ran(F)=\ran(M_w)^\perp$. 
In this case $F=\lambda -V$ where $V$ is a sum of a unilateral shift, a backward shift and truncated shifts of some indices ($V$ is a completely non-unitary power partial isometry see \cite{HW}). Indeed, let $\{{\bf j}_l\}_{l=l_0}^{l_1}$ be the increasing sequence of all $j$ such that $M_j>M_{j+1}$, where $-\infty\leq l_0\leq l_1\leq\infty.$ The sequence is nonempty, as the diagram is non-simple. Then, \begin{itemize}
     \item if $j_0=-\infty$ and $l_0\in\mathbb{Z}$, then $E_{b}:=\bigoplus_{j\leq {\bf j}_{l_0}}\mathbb{C}f_j$ reduces $V$ to a backward shift, denote $V_b=V|_{E_{b}},$
    \item if $j_1=\infty$ and $l_1\in\mathbb{Z}$, then $E_{s}:=\bigoplus_{j\geq {\bf j}_{l_1}}\mathbb{C}f_j$ reduces $V$ to a unilateral shift, denote $V_s=V|_{E_{s}},$
     \item if $l_1>l_0$, then $\mathcal{E}_l:=\bigoplus_{j={\bf j}_l+1}^{{\bf j}_{l+1}}\mathbb{C}f_j$ is finite dimensional space reducing $V$ to a truncated shift of index ${\bf j}_{l+1}-{\bf j}_l,$ for any $l_0-1< l<l_1,$ denote $V_k=V|_{E_{k}},$ where $E_k:=\bigoplus_{l:{\bf j}_{l+1}-{\bf j}_l=k}\mathcal{E}_l$. 
\end{itemize}  

 Note that $V_k$ is a truncated shift of multiplicity equal to cardinality of $\{l:{\bf j}_{l+1}-{\bf j}_l=k\}.$ 
 We may decompose $F=(\lambda-V_s)\oplus(\lambda-V_b)\oplus (\lambda -\bigoplus V_k)$ (see \cite{HW}). Note that $\ran(F)$ is dense and it is closed if and only if  $E_s=E_b=\{0\}$ and by Remark \ref{ts} if there is $K$ such that $E_k=\{0\}$ for any $k\ge K$.

We make the proof of this part for types of diagrams as in Case I. However, we get the joint result: $\ran(F)$ is not closed if and only if $\delta_-= 0$ or $\delta_+= 0.$

{\bf Case II.1}: $j_0=-\infty$, $j_1=\infty$. 

It is clear that if $l_0\in\mathbb{Z}$ (or $l_1\in\mathbb{Z}$) then $\delta_-=0$ and $E_b\neq\{0\}$ (resp. $\delta_+=0$ and $E_s\neq\{0\}$) which yields $\ran(F)$ is not closed. Hence we my restrict to the case $l_0=-\infty$ and $l_1=\infty.$ In such case, by Remark \ref{ts} we need to check that $\sup\{k:E_k\neq\{0\}\}=\infty$ if and only if  $\delta_-=0 $ or $\delta_+= 0.$  

Assume $\sup\{k:E_k\neq\{0\}\}=\infty$  and consider a sequence $\{k_i\}_{i\ge 0}$ such that $k_i\to\infty$ and $E_{k_i}\not=\{0\}$. Note that $E_{k_i}\not=\{0\}$ if and only if there is $l_{k_i}$ such that $\mathcal{E}_{l_{k_i}}\neq\{0\}$ and ${\bf j}_{l_{k_i}+1}-{\bf j}_{l_{k_i}}={k_i}$, so $M_{{\bf j}_{l_{k_i}}}>M_{{\bf j}_{l_{k_i}}+1}=M_{{\bf j}_{l_{k_i}}+2}=\ldots=M_{{\bf j}_{l_{k_i}+1}}>M_{{\bf j}_{l_{k_i}+1}+1}.$  Hence
$\frac{M_{{\bf j}_{l_{k_i}}+1}-M_{{\bf j}_{l_{k_i}}+n+1}}{n}=0$, for any $n\leq {k_i}-1$. Since for at most one $k_i$ may be ${\bf j}_{l_{k_i}}<0\leq {\bf j}_{l_{k_i}+1}$ we may choose a subsequence of $\{k_i\}$ such that all corresponding ${\bf j}_{l_{k_i}}, {\bf j}_{l_{k_i}+1}$ are negative or all are positive.
Hence \begin{itemize}
    \item for each $n$ there is $j\geq 0$ such that $M_{j+1}=M_{j+2}=\ldots=M_{j+n}$ or
    \item for each $n$ there is $j<-n$ such that $M_{j+1}=M_{j+2}=\ldots=M_{j+n}$.
\end{itemize}
Thus either $\inf_{j\ge 0} \frac{M_{j+1}-M_{n+j}}{n-1}=0$ for any $n>0,$ so $\delta_+=0$ or $\inf_{j\le 0} \frac{M_{j-1-n}-M_{j-1}}{n}=0$ for any $n>0$ so $\delta_-=0$.

For the reverse implication we assume $\delta_+=0$ (the case $\delta_-=0$ goes the same way).  Hence for any $K>0$ and arbitrary large $n$ say $n>K$ there is $j\geq-1$ such that $\frac{M_{j}-M_{n+j}}{n}<\frac1K$. Let us show that $E_k\neq\{0\}$ for some $k\ge\frac{K-1}{2}.$

Note that $M_{j}-M_{n+j}=\sum_{i=1}^{n}(M_{i+j-1}-M_{i+j})$ and the sequence $M_{j}-M_{j+1}, M_{j+1}-M_{j+2},\dots, M_{n+j-1}-M_{n+j}$ consist of zeros and positive integers. Let $m$ be the number of positive integers in the sequence. Since $m\leq \sum_{i=1}^{n}(M_{i+j-1}-M_{i+j})=M_{j}-M_{n+j}<\frac{n}{K}$ the number of zeros $n-m>n-\frac{n}{K}.$ Consequently there is a fragment of the sequence $M_{j}-M_{j+1}, M_{j+1}-M_{j+2},\dots, M_{n+j-1}-M_{n+j}$ consisting of zeros and of length not less than $\frac{n-m}{m+1}>\frac{n-\frac{n}{K}}{1+\frac{n}{K}}=\frac{K-1}{\frac{K}{n}+1}>\frac{K-1}{2}.$ Such fragment may be described by integers $\iota, \kappa$ such that $j\leq\iota<\iota+\kappa<n+j$ and $\kappa>\frac{K-1}{2}$ and $M_{\iota}-M_{\iota+1}=M_{\iota+1}-M_{\iota+2}=\dots=M_{\iota+\kappa-1}-M_{\iota+\kappa}=0,$ so $M_{\iota}=\dots=M_{\iota+\kappa}$. By the definition of $\{{\bf j}_l\}_{l=l_0}^{l_1}$ and $l_0=-\infty$, $l_1=\infty$ there is $l$ such that ${\bf j}_{l+1}=\max\{j:M_j=M_\iota\}$ and ${\bf j}_{l}=\min\{j:M_j=M_\iota\}-1.$ Hence $\{0\}\neq\mathcal{E}_l\subset E_k,$ where $k:={\bf j}_{l+1}-{\bf j}_{l}\ge\kappa> \frac{K-1}{2}.$  Since $K$ may be arbitrary large we get $\sup\{k:E_k\neq\{0\}\}=\infty.$

Summing up, $(0,\lambda)\in\Gamma_3(M_w,M_z)$ if and only if $\ran(F)$ is not closed if and only if  $\delta_-=0$ or $\delta_+=0$.

{\bf Case II.2} $j_1\in\ZZ, j_0=-\infty$. 

Note that, by definition $\delta_+=\infty.$
Similarly as in Case II.1 we can see that $(0,\lambda)\in\Gamma_3(M_w,M_z)$ if and only if $\delta_-=0$.

{\bf Case II.3} $j_0\in\ZZ, j_1=\infty$. 

By definition $\delta_-=\infty$ and similarly as in Case II.1 we can see that $(0,\lambda)\in\Gamma_3(M_w,M_z)$ if and only if $\delta_+=0$.

{\bf Case II.4} $j_0,j_1\in\ZZ$.

Then, $\delta_-=\delta_+=\infty$ and $F$ is an operator on finite dimensional space and $f_j\in\ran(F)$, for $j=j_0,j_0+1,\ldots,j_1$. Hence $\ran(F)=span\{f_j:j=j_0,j_0+1,\ldots,j_1\}=\ran(M_w)^\perp$. So $(\mu,\lambda)\not\in\Gamma_3(M_w,M_z)$.

All {\bf Case II} may be sum up with  $\{0\}\times\mathbb{T}\subset\Gamma_3(M_w,M_z)$ if and only if $\delta_-=0$ or $\delta_+=0.$ Otherwise $(\{0\}\times\mathbb{T})\cap\Gamma_3(M_w,M_z)=\emptyset.$

{\bf Step II}.

For the symmetric cases we take advantage of \begin{equation*}(\mu,\lambda)\in\Gamma_3(M_w,M_z) \iff (\lambda,\mu)\in\Gamma_3(M_z,M_w).\end{equation*} Moreover, $(M_w,M_z)$ on $H_J$ and  $(M_z,M_w)$ on $ H_{J'}$, for 
$J'=\{(j,i):(i,j)\in J\}$ are the same two operators (different by order as pairs) - see the pictures in Final step. Moreover, 
if $\delta_-,\delta_+,\rho_-,\rho_+$ are calculated  by \eqref{deltarho+-} for $J$ and $\delta'_-,\delta'_+,\rho'_-, \rho'_+$ are calculated  for $J'$, then
$\delta'_+=\frac1\rho_-$, $\rho'_+=\frac1\delta_-$, $\delta'_-=\frac1\rho_+$ and $\rho'_-=\frac1\delta_+$, where $\frac10=\infty$ (as the values are only non-negative) and $\frac1\infty=0.$ Indeed, it is enough to interpret any of the values in \eqref{deltarho+-} as a coefficient of a half-line. Such half-lines are uniquely determined by a diagram and so, they transforms via the same symmetry as $J$ transforms to $J'$.

{\bf Case I'}: $\mu\in\overline{\mathbb{D}}, 0<|\lambda|<1.$

Let us write down the result following from Case I applied to the pair $(M_z, M_w)$ on the space $H_{J'}$. We denote by $j'_0, j'_1$ values corresponding to $j_0, j_1$ calculated for the diagram $J'$.

{\bf Case I'.1}: $j'_0=-\infty, j'_1=\infty.$
We get $(\mu,\lambda)\in\Gamma_3(M_w, M_z)$ if and only if $$|\lambda|^{\frac1{\delta_+}}\leq|\mu|\leq |\lambda|^{\frac1{\rho_+}} \textnormal{ or } |\lambda|^{\frac1{\delta_-}}\leq|\mu|\leq |\lambda|^{\frac1{\rho_-}} \textnormal{ or } |\lambda|^{\frac1{\rho_+}}\leq|\mu|\leq |\lambda|^{\frac1{\delta_-}}.$$

{\bf Case I'.2}: $j'_0=-\infty, j'_1\in\mathbb{Z}.$

We get $(\mu,\lambda)\in\Gamma_3(M_w, M_z)$ if and only if $|\lambda|^{\frac1{\delta_+}}\leq|\mu|\leq |\lambda|^{\frac1{\rho_+}}.$

{\bf Case I'.3}: $j'_0\in\mathbb{Z}, j'_1=\infty.$

We get $(\mu,\lambda)\in\Gamma_3(M_w, M_z)$ if and only if $|\mu|\leq|\lambda|^\frac{1}{\rho_-}.$

{\bf Case I'.4}: $j'_0, j'_1\in\mathbb{Z}.$

We get $(\mu,\lambda)\in\Gamma_3(M_w, M_z)$ if and only if $\mu=0$.

{\bf Case II'}: $|\mu|=1$ and $\lambda=0.$

By Case II for the pair $(M_z, M_w)$ on the space $H_{J'}$ we get $\mathbb{T}\times\{0\}\subset\Gamma_3(M_w,M_z)$ if and only if $\delta'_-=\frac{1}{\rho_+}=0$ or $\delta'_+=\frac{1}{\rho_-}=0,$ so $\rho_-=\infty$ or $\rho_+=\infty.$ Otherwise, $(\mathbb{T}\times\{0\})\cap\Gamma_3(M_w,M_z)=\emptyset.$

{\bf Final step}

The sum of sets calculated in Cases I, II, I', II' and $\mathbb{T}^2\cup\{(0,0)\}$ is equal to  $\Gamma_3(M_w,M_z)\cup\mathbb{T}^2\cup\{(0,0)\}.$ However, to calculate such sum we need to combine types of diagrams from Case I, II, with types for Case I', II'. In other words, we need to consider all possible combinations - finite or $\pm\infty$ - of $j_0,j_1, j'_0, j'_1$. Note that: if $j_1'\in\mathbb{Z},$ then $j_0=-\infty$, if $j_0'\in\mathbb{Z},$ then $j_1=\infty,$ if $j_1\in\mathbb{Z},$ then $j'_0=-\infty$ and if $j_0\in\mathbb{Z},$ then $j'_1=\infty$. On the other hand, $j'_1=\infty$ does not say anything neither about $j_1$ nor $j_0$ and similarly other infinite values. The examples are in pictures below.

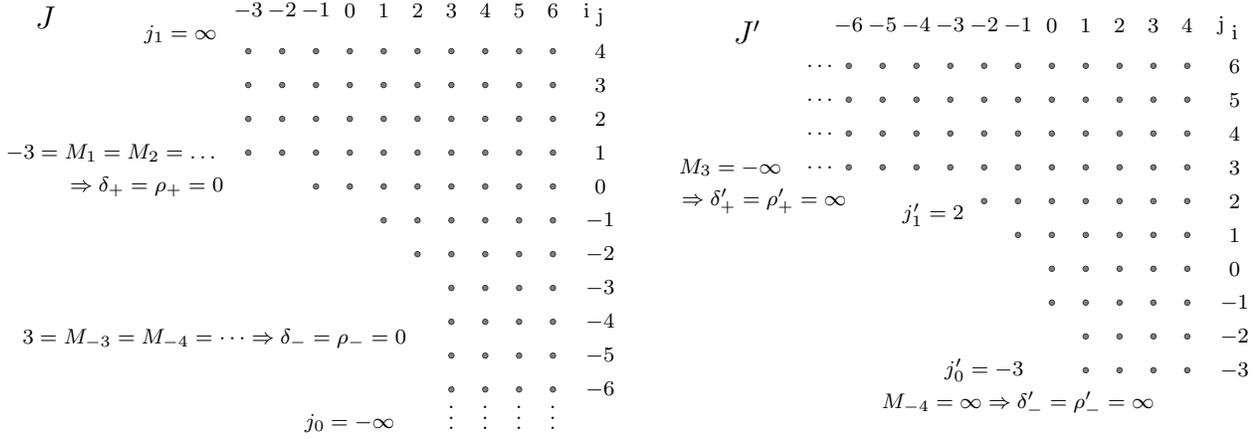
\begin{minipage}{0.45\linewidth}
\begin{center}
 \begin{tikzpicture}
[yscale=0.45, xscale=0.45, auto,
kropka/.style={draw=black!75, circle, inner sep=0pt, minimum size=2pt, fill=black!50},
kropa/.style={draw=black, circle, inner sep=0pt, minimum size=4pt, fill=black!100},
kolo/.style={draw=black, circle, inner sep=0pt, minimum size=8pt, semithick},
kwadrat/.style={draw=black, rectangle, inner sep=0pt,  minimum size=10pt, semithick},
romb/.style={draw=black, diamond, inner sep=0pt, minimum size=8pt, semithick},
trojkat/.style ={draw=black, regular polygon, regular polygon sides=3, inner sep=0pt, minimum size=8pt, semithick}]

\path (6,-4) node [kropka] {} (5,-4) node [kropka] {} (4,-4) node [kropka] {} (3,-4) node [kropka] {};
\path (6,-3) node [kropka] {} (5,-3) node [kropka] {} (4,-3) node [kropka] {} (3,-3) node [kropka] {};
\path (6,-2) node [kropka] {} (5,-2) node [kropka] {} (4,-2) node [kropka] {} (3,-2) node [kropka] {};
\path (6,-1) node [kropka] {} (5,-1) node [kropka] {} (4,-1) node [kropka] {} (3,-1) node [kropka] {};
\path (6,0) node [kropka] {} (5,0) node [kropka] {} (4,0) node [kropka] {} (3,0) node [kropka] {} (2,0) node [kropka] {};
\path (6,1) node [kropka] {} (5,1) node [kropka] {} (4,1) node [kropka] {} (3,1) node [kropka] {} (2,1) node [kropka] {}
(1,1) node [kropka] {};
\path (6,2) node [kropka] {} (5,2) node [kropka] {} (4,2) node [kropka] {} (3,2) node [kropka] {} (2,2) node [kropka] {}
(1,2) node [kropka] {} (0,2) node [kropka] {} (-1,2) node [kropka] {};
\path (6,3) node [kropka] {} (5,3) node [kropka] {} (4,3) node [kropka] {} (3,3) node [kropka] {} (2,3) node [kropka] {}
(1,3) node [kropka] {} (0,3) node [kropka] {} (-1,3) node [kropka] {} (-2,3) node [kropka] {} (-3,3) node [kropka] {};
\path (6,4) node [kropka] {} (5,4) node [kropka] {} (4,4) node [kropka] {} (3,4) node [kropka] {} (2,4) node [kropka] {}
(1,4) node [kropka] {} (0,4) node [kropka] {} (-1,4) node [kropka] {} (-2,4) node [kropka] {} (-3,4) node [kropka] {};
\path (6,5) node [kropka] {} (5,5) node [kropka] {} (4,5) node [kropka] {} (3,5) node [kropka] {} (2,5) node [kropka] {}
(1,5) node [kropka] {} (0,5) node [kropka] {} (-1,5) node [kropka] {} (-2,5) node [kropka] {} (-3,5) node [kropka] {};
\path (6,6) node [kropka] {} (5,6) node [kropka] {} (4,6) node [kropka] {} (3,6) node [kropka] {} (2,6) node [kropka] {}
(1,6) node [kropka] {} (0,6) node [kropka] {} (-1,6) node [kropka] {} (-2,6) node [kropka] {} (-3,6) node [kropka] {};

\node at  (3,-4.6) {\tiny{$\vdots$}};
\node at  (4,-4.6) {\tiny{$\vdots$}};
\node at  (5,-4.6) {\tiny{$\vdots$}};
\node at  (6,-4.6) {\tiny{$\vdots$}};
\node at (-3,7.2) {\tiny{$-3$}};
\node at (-2,7.2) {\tiny{$-2$}};
\node at (-1,7.2) {\tiny{$-1$}};
\node at (0,7.2) {\tiny{$0$}};
\node at (1,7.2) {\tiny{$1$}};
\node at (2,7.2) {\tiny{$2$}};
\node at (3,7.2) {\tiny{$3$}};
\node at (4,7.2) {\tiny{$4$}};
\node at (5,7.2) {\tiny{$5$}};
\node at (6,7.2) {\tiny{$6$}};
\node at (7,7.2) {\tiny{i}};

\node at (7.4, 7) {\tiny{j}};
\node at (7.4, 6) {\tiny{$4$}};
\node at (7.4, 5) {\tiny{$3$}};
\node at (7.4, 4) {\tiny{$2$}};
\node at (7.4, 3) {\tiny{$1$}};
\node at (7.4, 2) {\tiny{$0$}};
\node at (7.4, 1) {\tiny{$-1$}};
\node at (7.4, 0) {\tiny{$-2$}};
\node at (7.4, -1) {\tiny{$-3$}};
\node at (7.4, -2) {\tiny{$-4$}};
\node at (7.4, -3) {\tiny{$-5$}};
\node at (7.4, -4) {\tiny{$-6$}};

\node at (0, -5) {\tiny{$j_0=-\infty$}};
\node at (-5, 6.5) {\tiny{$j_1=\infty$}};
\node at (-7,3) {\tiny{$-3=M_1=M_2=\dots$}};
\node at (-6,2) {\tiny{$\Rightarrow \delta_+=\rho_+=0$}};
\node at (-4,-2.5) {\tiny{$3=M_{-3}=M_{-4}=\dots\Rightarrow \delta_-=\rho_-=0$}};
\node at (-9, 7) {$J$};
\end{tikzpicture}
\end{center}
\end{minipage}
\begin{minipage}{0.05\linewidth}
\end{minipage}
\begin{minipage}{0.5\linewidth}
\begin{center}
 \begin{tikzpicture}
[yscale=0.45, xscale=0.45, auto,
kropka/.style={draw=black!75, circle, inner sep=0pt, minimum size=2pt, fill=black!50},
kropa/.style={draw=black, circle, inner sep=0pt, minimum size=4pt, fill=black!100},
kolo/.style={draw=black, circle, inner sep=0pt, minimum size=8pt, semithick},
kwadrat/.style={draw=black, rectangle, inner sep=0pt,  minimum size=10pt, semithick},
romb/.style={draw=black, diamond, inner sep=0pt, minimum size=8pt, semithick},
trojkat/.style ={draw=black, regular polygon, regular polygon sides=3, inner sep=0pt, minimum size=8pt, semithick}]

\path (6,-3) node [kropka] {} (5,-3) node [kropka] {} (4,-3) node [kropka] {} (3,-3) node [kropka] {};
\path (6,-2) node [kropka] {} (5,-2) node [kropka] {} (4,-2) node [kropka] {} (3,-2) node [kropka] {};
\path (6,-1) node [kropka] {} (5,-1) node [kropka] {} (4,-1) node [kropka] {} (3,-1) node [kropka] {} (2,-1) node [kropka] {};
\path (6,0) node [kropka] {} (5,0) node [kropka] {} (4,0) node [kropka] {} (3,0) node [kropka] {} (2,0) node [kropka] {} ;
\path (6,1) node [kropka] {} (5,1) node [kropka] {} (4,1) node [kropka] {} (3,1) node [kropka] {} (2,1) node [kropka] {}
(1,1) node [kropka] {};
\path (6,2) node [kropka] {} (5,2) node [kropka] {} (4,2) node [kropka] {} (3,2) node [kropka] {} (2,2) node [kropka] {}
(1,2) node [kropka] {} (0,2) node [kropka] {};
\path (6,3) node [kropka] {} (5,3) node [kropka] {} (4,3) node [kropka] {} (3,3) node [kropka] {} (2,3) node [kropka] {}
(1,3) node [kropka] {} (0,3) node [kropka] {} (-1,3) node [kropka] {} (-2,3) node [kropka] {} (-3,3) node [kropka] {} (-4,3) node [kropka] {};
\path (6,4) node [kropka] {} (5,4) node [kropka] {} (4,4) node [kropka] {} (3,4) node [kropka] {} (2,4) node [kropka] {}
(1,4) node [kropka] {} (0,4) node [kropka] {} (-1,4) node [kropka] {} (-2,4) node [kropka] {} (-3,4) node [kropka] {} (-4,4) node [kropka] {};
\path (6,5) node [kropka] {} (5,5) node [kropka] {} (4,5) node [kropka] {} (3,5) node [kropka] {} (2,5) node [kropka] {}
(1,5) node [kropka] {} (0,5) node [kropka] {} (-1,5) node [kropka] {} (-2,5) node [kropka] {} (-3,5) node [kropka] {} (-4,5) node [kropka] {};
\path (6,6) node [kropka] {} (5,6) node [kropka] {} (4,6) node [kropka] {} (3,6) node [kropka] {} (2,6) node [kropka] {}
(1,6) node [kropka] {} (0,6) node [kropka] {} (-1,6) node [kropka] {} (-2,6) node [kropka] {} (-3,6) node [kropka] {} (-4,6) node [kropka] {};

\node at  (-4.8,6) {\tiny{$\dots$}};
\node at  (-4.8,5) {\tiny{$\dots$}};
\node at  (-4.8,4) {\tiny{$\dots$}};
\node at  (-4.8,3) {\tiny{$\dots$}};
\node at (-4,7.2) {\tiny{$-6$}};
\node at (-3,7.2) {\tiny{$-5$}};
\node at (-2,7.2) {\tiny{$-4$}};
\node at (-1,7.2) {\tiny{$-3$}};
\node at (0,7.2) {\tiny{$-2$}};
\node at (1,7.2) {\tiny{$-1$}};
\node at (2,7.2) {\tiny{$0$}};
\node at (3,7.2) {\tiny{$1$}};
\node at (4,7.2) {\tiny{$2$}};
\node at (5,7.2) {\tiny{$3$}};
\node at (6,7.2) {\tiny{$4$}};
\node at (7,7.2) {\tiny{j}};

\node at (7.4, 7) {\tiny{i}};
\node at (7.4, 6) {\tiny{$6$}};
\node at (7.4, 5) {\tiny{$5$}};
\node at (7.4, 4) {\tiny{$4$}};
\node at (7.4, 3) {\tiny{$3$}};
\node at (7.4, 2) {\tiny{$2$}};
\node at (7.4, 1) {\tiny{$1$}};
\node at (7.4, 0) {\tiny{$0$}};
\node at (7.4, -1) {\tiny{$-1$}};
\node at (7.4, -2) {\tiny{$-2$}};
\node at (7.4, -3) {\tiny{$-3$}};

\node at (0, -3) {\tiny{$j'_0=-3$}};
\node at (-1.5, 1.6) {\tiny{$j'_1=2$}};
\node at (-7, 7) {$J'$};
\node at (-7.5, 3) {\tiny{$M_3=-\infty$}};
\node at (-6.5, 2) {\tiny{$\Rightarrow \delta'_+=\rho'_+=\infty$}};
\node at (1, -4) {\tiny{$M_{-4}=\infty\Rightarrow \delta'_-=\rho'_-=\infty$}};
\end{tikzpicture}
\end{center}
\end{minipage}

In the following picture we assume the respective part of the diagram lies on or on the right of the dashed line. Hence there follows values of  $\delta_+,\rho_+,\delta'_-,\rho'_-.$

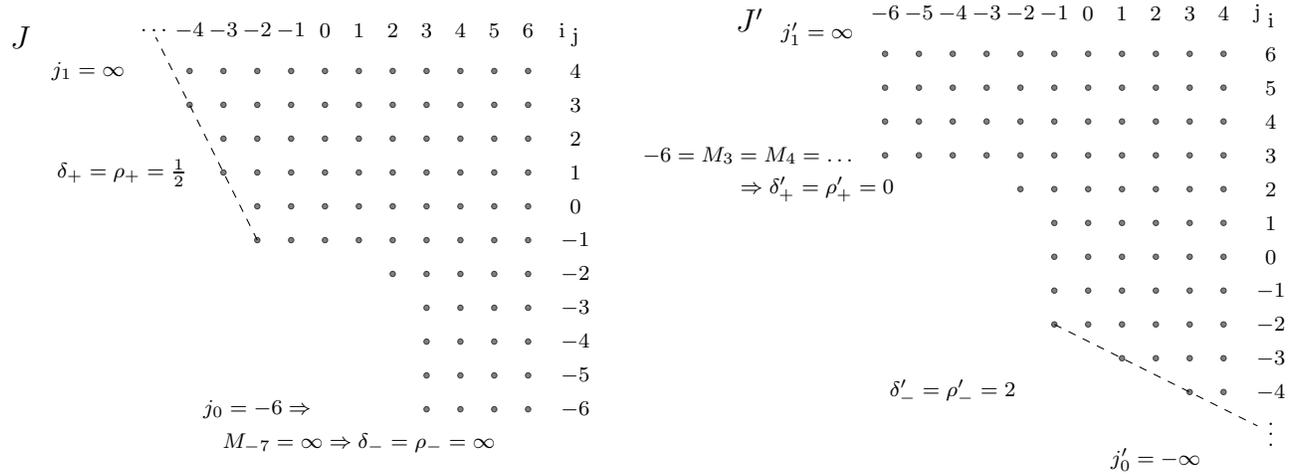
\begin{minipage}{0.45\linewidth}
\begin{center}
 \begin{tikzpicture}
[yscale=0.45, xscale=0.45, auto,
kropka/.style={draw=black!75, circle, inner sep=0pt, minimum size=2pt, fill=black!50},
kropa/.style={draw=black, circle, inner sep=0pt, minimum size=4pt, fill=black!100},
kolo/.style={draw=black, circle, inner sep=0pt, minimum size=8pt, semithick},
kwadrat/.style={draw=black, rectangle, inner sep=0pt,  minimum size=10pt, semithick},
romb/.style={draw=black, diamond, inner sep=0pt, minimum size=8pt, semithick},
trojkat/.style ={draw=black, regular polygon, regular polygon sides=3, inner sep=0pt, minimum size=8pt, semithick}]

\path (6,-4) node [kropka] {} (5,-4) node [kropka] {} (4,-4) node [kropka] {} (3,-4) node [kropka] {};
\path (6,-3) node [kropka] {} (5,-3) node [kropka] {} (4,-3) node [kropka] {} (3,-3) node [kropka] {};
\path (6,-2) node [kropka] {} (5,-2) node [kropka] {} (4,-2) node [kropka] {} (3,-2) node [kropka] {};
\path (6,-1) node [kropka] {} (5,-1) node [kropka] {} (4,-1) node [kropka] {} (3,-1) node [kropka] {};
\path (6,0) node [kropka] {} (5,0) node [kropka] {} (4,0) node [kropka] {} (3,0) node [kropka] {} (2,0) node [kropka] {};
\path (6,1) node [kropka] {} (5,1) node [kropka] {} (4,1) node [kropka] {} (3,1) node [kropka] {} (2,1) node [kropka] {}
(1,1) node [kropka] {} (0,1) node [kropka] {} (-1,1) node [kropka] {} (-2,1) node [kropka] {};
\path (6,2) node [kropka] {} (5,2) node [kropka] {} (4,2) node [kropka] {} (3,2) node [kropka] {} (2,2) node [kropka] {}
(1,2) node [kropka] {} (0,2) node [kropka] {} (-1,2) node [kropka] {} (-2,2) node [kropka] {};
\path (6,3) node [kropka] {} (5,3) node [kropka] {} (4,3) node [kropka] {} (3,3) node [kropka] {} (2,3) node [kropka] {}
(1,3) node [kropka] {} (0,3) node [kropka] {} (-1,3) node [kropka] {} (-2,3) node [kropka] {} (-3,3) node [kropka] {};
\path (6,4) node [kropka] {} (5,4) node [kropka] {} (4,4) node [kropka] {} (3,4) node [kropka] {} (2,4) node [kropka] {}
(1,4) node [kropka] {} (0,4) node [kropka] {} (-1,4) node [kropka] {} (-2,4) node [kropka] {} (-3,4) node [kropka] {};
\path (6,5) node [kropka] {} (5,5) node [kropka] {} (4,5) node [kropka] {} (3,5) node [kropka] {} (2,5) node [kropka] {}
(1,5) node [kropka] {} (0,5) node [kropka] {} (-1,5) node [kropka] {} (-2,5) node [kropka] {} (-3,5) node [kropka] {}(-4,5) node [kropka] {};
\path (6,6) node [kropka] {} (5,6) node [kropka] {} (4,6) node [kropka] {} (3,6) node [kropka] {} (2,6) node [kropka] {}
(1,6) node [kropka] {} (0,6) node [kropka] {} (-1,6) node [kropka] {} (-2,6) node [kropka] {} (-3,6) node [kropka] {} (-4,6) node [kropka] {};
\draw [dashed] (-2,1)--(-5,7);

\node at (-5,7.2) {\tiny{$\dots$}};
\node at (-4,7.2) {\tiny{$-4$}};
\node at (-3,7.2) {\tiny{$-3$}};
\node at (-2,7.2) {\tiny{$-2$}};
\node at (-1,7.2) {\tiny{$-1$}};
\node at (0,7.2) {\tiny{$0$}};
\node at (1,7.2) {\tiny{$1$}};
\node at (2,7.2) {\tiny{$2$}};
\node at (3,7.2) {\tiny{$3$}};
\node at (4,7.2) {\tiny{$4$}};
\node at (5,7.2) {\tiny{$5$}};
\node at (6,7.2) {\tiny{$6$}};
\node at (7,7.2) {\tiny{i}};

\node at (7.4, 7) {\tiny{j}};
\node at (7.4, 6) {\tiny{$4$}};
\node at (7.4, 5) {\tiny{$3$}};
\node at (7.4, 4) {\tiny{$2$}};
\node at (7.4, 3) {\tiny{$1$}};
\node at (7.4, 2) {\tiny{$0$}};
\node at (7.4, 1) {\tiny{$-1$}};
\node at (7.4, 0) {\tiny{$-2$}};
\node at (7.4, -1) {\tiny{$-3$}};
\node at (7.4, -2) {\tiny{$-4$}};
\node at (7.4, -3) {\tiny{$-5$}};
\node at (7.4, -4) {\tiny{$-6$}};

\node at (-7, 6) {\tiny{$j_1=\infty$}};
\node at (-6,3) {\tiny{$\delta_+=\rho_+=\frac12$}};
\node at (-2,-4) {\tiny{$j_0=-6\Rightarrow $}};
\node at (1,-5) {\tiny{$M_{-7}=\infty\Rightarrow \delta_-=\rho_-=\infty$}};
\node at (-9, 7) {$J$};
\end{tikzpicture}
\end{center}
\end{minipage}
\begin{minipage}{0.05\linewidth}
\end{minipage}
\begin{minipage}{0.5\linewidth}
\begin{center}
 \begin{tikzpicture}
[yscale=0.45, xscale=0.45, auto,
kropka/.style={draw=black!75, circle, inner sep=0pt, minimum size=2pt, fill=black!50},
kropa/.style={draw=black, circle, inner sep=0pt, minimum size=4pt, fill=black!100},
kolo/.style={draw=black, circle, inner sep=0pt, minimum size=8pt, semithick},
kwadrat/.style={draw=black, rectangle, inner sep=0pt,  minimum size=10pt, semithick},
romb/.style={draw=black, diamond, inner sep=0pt, minimum size=8pt, semithick},
trojkat/.style ={draw=black, regular polygon, regular polygon sides=3, inner sep=0pt, minimum size=8pt, semithick}]

\path (6,-4) node [kropka] {} (5,-4) node [kropka] {} ;
\path (6,-3) node [kropka] {} (5,-3) node [kropka] {} (4,-3) node [kropka] {} (3,-3) node [kropka] {};
\path (6,-2) node [kropka] {} (5,-2) node [kropka] {} (4,-2) node [kropka] {} (3,-2) node [kropka] {}(2,-2) node [kropka] {}(1,-2) node [kropka] {};
\path (6,-1) node [kropka] {} (5,-1) node [kropka] {} (4,-1) node [kropka] {} (3,-1) node [kropka] {} (2,-1) node [kropka] {}(1,-1) node [kropka] {} ;
\path (6,0) node [kropka] {} (5,0) node [kropka] {} (4,0) node [kropka] {} (3,0) node [kropka] {} (2,0) node [kropka] {} (1,0) node [kropka] {};
\path (6,1) node [kropka] {} (5,1) node [kropka] {} (4,1) node [kropka] {} (3,1) node [kropka] {} (2,1) node [kropka] {}
(1,1) node [kropka] {};
\path (6,2) node [kropka] {} (5,2) node [kropka] {} (4,2) node [kropka] {} (3,2) node [kropka] {} (2,2) node [kropka] {}
(1,2) node [kropka] {} (0,2) node [kropka] {};
\path (6,3) node [kropka] {} (5,3) node [kropka] {} (4,3) node [kropka] {} (3,3) node [kropka] {} (2,3) node [kropka] {}
(1,3) node [kropka] {} (0,3) node [kropka] {} (-1,3) node [kropka] {} (-2,3) node [kropka] {} (-3,3) node [kropka] {} (-4,3) node [kropka] {};
\path (6,4) node [kropka] {} (5,4) node [kropka] {} (4,4) node [kropka] {} (3,4) node [kropka] {} (2,4) node [kropka] {}
(1,4) node [kropka] {} (0,4) node [kropka] {} (-1,4) node [kropka] {} (-2,4) node [kropka] {} (-3,4) node [kropka] {} (-4,4) node [kropka] {};
\path (6,5) node [kropka] {} (5,5) node [kropka] {} (4,5) node [kropka] {} (3,5) node [kropka] {} (2,5) node [kropka] {}
(1,5) node [kropka] {} (0,5) node [kropka] {} (-1,5) node [kropka] {} (-2,5) node [kropka] {} (-3,5) node [kropka] {} (-4,5) node [kropka] {};
\path (6,6) node [kropka] {} (5,6) node [kropka] {} (4,6) node [kropka] {} (3,6) node [kropka] {} (2,6) node [kropka] {}
(1,6) node [kropka] {} (0,6) node [kropka] {} (-1,6) node [kropka] {} (-2,6) node [kropka] {} (-3,6) node [kropka] {} (-4,6) node [kropka] {};
\draw [dashed] (1,-2)--(7,-5);

\node at (-4,7.2) {\tiny{$-6$}};
\node at (-3,7.2) {\tiny{$-5$}};
\node at (-2,7.2) {\tiny{$-4$}};
\node at (-1,7.2) {\tiny{$-3$}};
\node at (0,7.2) {\tiny{$-2$}};
\node at (1,7.2) {\tiny{$-1$}};
\node at (2,7.2) {\tiny{$0$}};
\node at (3,7.2) {\tiny{$1$}};
\node at (4,7.2) {\tiny{$2$}};
\node at (5,7.2) {\tiny{$3$}};
\node at (6,7.2) {\tiny{$4$}};
\node at (7,7.2) {\tiny{j}};

\node at (7.4, 7) {\tiny{i}};
\node at (7.4, 6) {\tiny{$6$}};
\node at (7.4, 5) {\tiny{$5$}};
\node at (7.4, 4) {\tiny{$4$}};
\node at (7.4, 3) {\tiny{$3$}};
\node at (7.4, 2) {\tiny{$2$}};
\node at (7.4, 1) {\tiny{$1$}};
\node at (7.4, 0) {\tiny{$0$}};
\node at (7.4, -1) {\tiny{$-1$}};
\node at (7.4, -2) {\tiny{$-2$}};
\node at (7.4, -3) {\tiny{$-3$}};
\node at (7.4, -4) {\tiny{$-4$}};
\node at (7.4, -5) {\tiny{$\vdots$}};
\node at (4, -6) {\tiny{$j'_0=-\infty$}};
\node at (-6, 6.6) {\tiny{$j'_1=\infty$}};
\node at (-8, 7) {$J'$};
\node at (-8, 3) {\tiny{$-6=M_3=M_4=\dots$}};
\node at (-6, 2) {\tiny{$\Rightarrow \delta'_+=\rho'_+=0$}};
\node at (-2, -4) {\tiny{$\delta'_-=\rho'_-=2$}};
\end{tikzpicture}
\end{center}
\end{minipage}

{\bf Case A}: $j_0=j'_0=-\infty, j_1=j_1'=\infty.$

Since $j_1=j'_1=\infty$ none of $M_w, M_z$ has unitary part, so we need to show:
$$\Gamma_3(M_w,M_z)\cup\mathbb{T}^2=\{(\mu,\lambda)\in\overline{\mathbb{D}}^2:
|\mu|^{\rho_-}\leq|\lambda|\leq |\mu|^{\delta_-}   \textnormal{ or }
|\mu|^{\delta_-}\leq|\lambda|\leq |\mu|^{\rho_+}  \textnormal{ or }
|\mu|^{\rho_+}\leq|\lambda|\leq |\mu|^{\delta_+}\}.$$

By equivalence of \eqref{one} with \eqref{two} in Remark \ref{pab} the sum of sets in Case I and Case I' (precisely Case I.1 and Case I'.1) may be represented as
$$\{(\mu,\lambda)\in(\mathbb{D}\times \overline{\mathbb{D}})\cup (\overline{\mathbb{D}}\times \mathbb{D}):
|\mu|^{\rho_-}\leq|\lambda|\leq |\mu|^{\delta_-}   \textnormal{ or }
|\mu|^{\delta_-}\leq|\lambda|\leq |\mu|^{\rho_+}  \textnormal{ or }
|\mu|^{\rho_+}\leq|\lambda|\leq |\mu|^{\delta_+}\}.$$
 By Remark \ref{pab} $(\mu,\lambda)\in\{0\}\times \mathbb{T}$ satisfy any pair of inequalities in the set above if and only if $\delta_-=0$ or $\delta_+=0$ which is precisely as in Case II. Indeed, $\rho_+\ge\delta_+\ge 0,$ so $\rho_+=0$ does not give an extra case. Hence $$\{(\mu,\lambda)\in(\mathbb{D}\times \overline{\mathbb{D}})\cup (\overline{\mathbb{D}}\times \mathbb{D})\cup \{0\}\times \mathbb{T}:
|\mu|^{\rho_-}\leq|\lambda|\leq |\mu|^{\delta_-}   \textnormal{ or }
|\mu|^{\delta_-}\leq|\lambda|\leq |\mu|^{\rho_+}  \textnormal{ or }
|\mu|^{\rho_+}\leq|\lambda|\leq |\mu|^{\delta_+}\}$$ is the sum of sets in Cases I, I' and II. We do the same with Case II'.  Consequently the sum of results of Case I, II, I', II' and $\{(0,0)\}\cup\mathbb{T}^2$ is indeed equal to
$$\{(\mu,\lambda)\in\overline{\mathbb{D}}^2:
|\mu|^{\rho_-}\leq|\lambda|\leq |\mu|^{\delta_-}   \textnormal{ or }
|\mu|^{\delta_-}\leq|\lambda|\leq |\mu|^{\rho_+}  \textnormal{ or }
|\mu|^{\rho_+}\leq|\lambda|\leq |\mu|^{\delta_+}\}.$$
 {\bf Case B}: $j_0\in\mathbb{Z}, j'_0=-\infty, j_1=j_1'=\infty.$

 Since $j_0\in\mathbb{Z}$ we get $\rho_-=\delta_-=\infty$. Moreover, $j_1=j'_1=\infty$ yields none of $M_w, M_z$ has unitary part, so as in Case A we need to show:
$$\Gamma_3(M_w,M_z)\cup\mathbb{T}^2=\{(\mu,\lambda)\in\overline{\mathbb{D}}^2:
|\mu|^{\rho_-}\leq|\lambda|\leq |\mu|^{\delta_-}   \textnormal{ or }
|\mu|^{\delta_-}\leq|\lambda|\leq |\mu|^{\rho_+}  \textnormal{ or }
|\mu|^{\rho_+}\leq|\lambda|\leq |\mu|^{\delta_+}\}.$$
 
Since $\rho_-=\delta_-=\infty$ the right hand side may be simplified. We use Remark \ref{pab}. The first condition $|\mu|^{\rho_-}\leq|\lambda|\leq |\mu|^{\delta_-}$   gives the set $\overline{\mathbb{D}}\times\{0\}\cup\mathbb{T}\times\overline{\mathbb{D}}.$ The second condition reduces to $|\lambda|\leq |\mu|^{\rho_+}$ which together with the third condition gives  $|\lambda|\leq |\mu|^{\delta_+}$ (as $\rho_+\geq\delta_+$). Since the last inequality is satisfied by values in $\overline{\mathbb{D}}\times\{0\}\cup\mathbb{T}\times\overline{\mathbb{D}}$ the whole set on the right may be described as 
$$\{(\mu,\lambda)\in\overline{\mathbb{D}}^2: |\lambda|\leq|\mu|^{\delta_+}\}.$$

Let us now calculate $\Gamma_3(M_w, M_z)\cup\mathbb{T}^2$. 
By Case I.3 we get the set $$\{(\mu,\lambda)\in\mathbb{D}\setminus\{0\}\times\overline{\mathbb{D}}:|\lambda|\leq |\mu|^{\delta_+}\}.$$
  By equivalence of \eqref{one} with \eqref{two} in Remark \ref{pab} the set in Case I'.1 is equal to
$$\{(\mu,\lambda)\in\overline{\mathbb{D}}\times \mathbb{D}\setminus\{0\}:|\mu|^{\rho_-}\leq|\lambda|\leq |\mu|^{\delta_-}   \textnormal{ or }
|\mu|^{\delta_-}\leq|\lambda|\leq |\mu|^{\rho_+}  \textnormal{ or }
|\mu|^{\rho_+}\leq|\lambda|\leq |\mu|^{\delta_+}\}$$
$$=\{(\mu,\lambda)\in\overline{\mathbb{D}}\times \mathbb{D}\setminus\{0\}:
|\lambda|\leq |\mu|^{\rho_+}  \textnormal{ or }
|\mu|^{\rho_+}\leq|\lambda|\leq |\mu|^{\delta_+}\}=\{(\mu,\lambda)\in\overline{\mathbb{D}}\times \mathbb{D}\setminus\{0\}:
|\lambda|\leq |\mu|^{\delta_+}\}.$$
Similarly like in Case A, sum of sets in Case I.3 and Case I'.1 is
$$\{(\mu,\lambda)\in\mathbb{D}\setminus\{0\}\times\overline{\mathbb{D}}\cup \overline{\mathbb{D}}\times \mathbb{D}\setminus\{0\}:|\lambda|\leq |\mu|^{\delta_+}\}.$$
Since $\rho_-=\infty$  Case II' gives the set $\mathbb{T}\times\{0\}.$ Note that $(\mu,\lambda)\in\mathbb{T}\times\{0\}$ satisfies the inequality $|\lambda|\leq|\mu|^{\delta_+}.$ 
Since $\delta_-=\infty\neq 0$ Case II gives the set $\{0\}\times\mathbb{T}$ if $\delta_+=0$ and $\emptyset$ otherwise. On the other hand,  $(\mu,\lambda)\in\{0\}\times\mathbb{T}$ satisfy  $|\lambda|\leq|\mu|^{\delta_+}$ if and only if $\delta_+=0.$ 

Hence the sum of sets calculated in Cases I, II, I', II' and $\{(0,0)\}\cup\mathbb{T}^2$ is indeed equal to
$$\{(\mu,\lambda)\in\overline{\mathbb{D}}^2:
|\lambda|\leq |\mu|^{\delta_+}\}.$$

  {\bf Case C}: $j_0, j_1\in\mathbb{Z}, j'_0=-\infty, j_1'=\infty.$

 Since $j_0, j_1\in\mathbb{Z}$ we get $\rho_+=\delta_+=\rho_-=\delta_-=\infty$ and $M_w$ has nontrivial Wold decomposition and $M_z$ is a unilateral shift. Hence we need to show $$\Gamma_3(M_w,M_z)\cup\mathbb{T}^2=\{(\mu,\lambda)\in\overline{\mathbb{D}}^2: |\mu|^{\rho_-}\leq |\lambda|\leq |\mu|^{\delta_-} \}\cup\mathbb{T}\times\cD.$$
Since $\rho_-=\delta_-=\infty$ the right hand side set may be simplified by Remark \ref{pab} to $$\mathbb{T}\times\overline{\mathbb{D}}\cup \mathbb{D}\times \{0\}.$$

Let us now calculate $\Gamma_3(M_w, M_z)\cup\mathbb{T}^2$. By Case I.4 we get $\mathbb{D}\setminus\{0\}\times \{0\}.$
By $\rho_+=\delta_+=\rho_-=\delta_-=\infty$ and Remark \ref{pab} the set calculated in Case I'.1 is equal to
$\mathbb{T}\times\mathbb{D}\setminus\{0\},$
 the set calculated in Case II is $\emptyset$ and the set calculated in Case II' is $\mathbb{T}\times\{0\}.$
The sum of the above sets and $\{(0,0)\}\cup\mathbb{T}^2$ is indeed equal to 
$$\mathbb{T}\times\overline{\mathbb{D}}\cup \mathbb{D}\times \{0\}.$$
{\bf Case D}: $j_0, j'_0\in\mathbb{Z}, j_1=j_1'=\infty.$

Since $j_0\in\mathbb{Z}$ we get $\rho_-=\delta_-=\infty$. Since $j'_0\in\mathbb{Z}$ we get $\frac{1}{\delta_+}=\rho'_-=\frac{1}{\rho_+}=\delta'_-=\infty,$ so $\delta_+=\rho_+=0$. Moreover, $j_1=j'_1=\infty$ yields none of $M_w, M_z$ has unitary part, so we need to show
$$\Gamma_3(M_w,M_z)\cup\mathbb{T}^2=\{(\mu,\lambda)\in\overline{\mathbb{D}}^2:
|\mu|^{\rho_-}\leq|\lambda|\leq |\mu|^{\delta_-}   \textnormal{ or }
|\mu|^{\delta_-}\leq|\lambda|\leq |\mu|^{\rho_+}  \textnormal{ or }
|\mu|^{\rho_+}\leq|\lambda|\leq |\mu|^{\delta_+}\}.$$

Note that by $\delta_-=\infty$ and $\rho_+=0$ and Remark \ref{pab} the middle term $|\mu|^{\delta_-}\leq|\lambda|\leq |\mu|^{\rho_+}$ above is satisfied on the whole $\overline{\mathbb{D}}^2.$ Hence we need to show 
$$\Gamma_3(M_w,M_z)\cup\mathbb{T}^2=\overline{\mathbb{D}}^2.$$
By Remark \ref{pab} inequalities in Case I.3 are satisfied on the whole domain $\mathbb{D}\setminus\{0\}\times\overline{\mathbb{D}}$ and similarly in the Case I'.3 where the domain is $\overline{\mathbb{D}}\times\mathbb{D}\setminus\{0\}.$ Sets calculated in Cases II and II' are $\{0\}\times\mathbb{T}$ and $\mathbb{T}\times\{0\}$ respectively.

The sum of all the sets above and $\{(0,0)\}\cup\mathbb{T}^2$ is indeed the whole $\overline{\mathbb{D}}^2.$

{\bf Case E}: $j_0=j'_0=-\infty, j_1=\infty , j_1'\in\mathbb{Z}.$

Since $j'_1\in\mathbb{Z}$ we get $\delta'_+=\frac{1}{\rho_-}=\rho'_+=\frac{1}{\delta_-}=\infty,$ so $\delta_-=\rho_-=0$ and $M_w$ is a unilateral shift and $M_z$ has nontrivial Wold decomposition.

Hence we need to show  $$\Gamma_3(M_w,M_z)\cup\mathbb{T}^2=\{(\mu,\lambda)\in \overline{\mathbb{D}}^2:|\mu|^{\rho_+}\leq|\lambda|\leq |\mu|^{\delta_+}\}\cup\cD\times\mathbb{T}.$$
Note that $\mathbb{T}\times\{0\}$ is a subset of the set above if and only if $\rho_+=\infty$.

By $\delta_-=\rho_-=0$ the set calculated in the Case I.1 is equal to
$$\{(\mu,\lambda)\in\mathbb{D}\setminus\{0\}\times \overline{\mathbb{D}}:|\lambda|=1  \textnormal{ or }
|\mu|^{\rho_+}\leq|\lambda|\leq |\mu|^{\delta_+}\}.$$
The set calculated in Case I'.2 is equal to $$\{(\mu,\lambda)\in \overline{\mathbb{D}}\times\mathbb{D}\setminus\{0\}:
|\mu|^{\rho_+}\leq|\lambda|\leq |\mu|^{\delta_+}\}.$$ 
Note that the alternate condition $|\lambda|=1$ is out of the domain of the latter set. Hence the sum of sets calculated in Cases I and I' is
$$\{(\mu,\lambda)\in(\mathbb{D}\setminus\{0\}\times \overline{\mathbb{D}})\cup(\overline{\mathbb{D}}\times\mathbb{D}\setminus\{0\}):|\lambda|=1  \textnormal{ or }
|\mu|^{\rho_+}\leq|\lambda|\leq |\mu|^{\delta_+}\}.$$
The set calculated in Case II is equal to $\{0\}\times\mathbb{T}\subset \Gamma_3(M_w,M_z)$ and by Remark \ref{pab} its elements satisfy $|\mu|^{\rho_+}\leq|\lambda|\leq |\mu|^{\delta_+}$. Hence it may be added to the previous sum by extending the domain by $\{0\}\times\mathbb{T},$ that is 
$$\{(\mu,\lambda)\in(\mathbb{D}\setminus\{0\}\times \overline{\mathbb{D}})\cup(\overline{\mathbb{D}}\times\mathbb{D}\setminus\{0\})\cup(\{0\}\times\mathbb{T}):|\lambda|=1  \textnormal{ or }
|\mu|^{\rho_+}\leq|\lambda|\leq |\mu|^{\delta_+}\}.$$
 The set calculated in  Case II' is equal to $\mathbb{T}\times\{0\}\subset \Gamma_3(M_w,M_z)$ if  $\rho_+=\infty$ and $\emptyset$ otherwise. On the other hand elements of $\mathbb{T}\times\{0\}$ satisfy $|\mu|^{\rho_+}\leq|\lambda|\leq |\mu|^{\delta_+}$ if and only if $\rho_+=\infty$. Therefore, again we may add this part to the previous sum by extending the domain by $\mathbb{T}\times\{0\}$. One can check that  the sum  above with $\{(0,0)\}\cup\mathbb{T}^2$ is indeed equal to 
$$\{(\mu,\lambda)\in \overline{\mathbb{D}}^2:|\mu|^{\rho_+}\leq|\lambda|\leq |\mu|^{\delta_+}   \textnormal{ or } |\lambda|=1\}.$$

{\bf Case F}: $j_0=j'_0=-\infty , j_1,j_1'\in\mathbb{Z}.$

In this case we get $\rho_-=\delta_-=0$ and $\rho_+=\delta_+=\infty.$ Moreover, both operators $M_w, M_z$ has nontrivial Wold decomposition. 

By the Case I.2 we get the set
$$\{(\mu,\lambda)\in\mathbb{D}\setminus\{0\}\times \overline{\mathbb{D}}:|\mu|^{\rho_-}\leq|\lambda|\leq |\mu|^{\delta_-}\}=
\mathbb{D}\setminus\{0\}\times\mathbb{T}.$$

By the Case I'.2 we get the set $\{(\mu,\lambda)\in \overline{\mathbb{D}}\times\mathbb{D}\setminus\{0\}:
|\lambda|^{\frac{1}{\delta_+}}\leq|\mu|\leq |\lambda|^{\frac{1}{\rho_+}}\},$ which by Remark \ref{pab} equals to  $$\{(\mu,\lambda)\in \overline{\mathbb{D}}\times\mathbb{D}\setminus\{0\}:
|\mu|^{\rho_+}\leq|\lambda|\leq |\mu|^{\delta_+}\}=\mathbb{T}\times\mathbb{D}\setminus\{0\}.$$ 

By Cases II and II' we get $\mathbb{T}\times\{0\}$ and $\{0\}\times\mathbb{T}$. 

Hence $$\Gamma_3(M_w,M_z)\cup\{(0,0\}\cup\mathbb{T}^2=\mathbb{T}\times\mathbb{D}\cup \mathbb{D}\times\mathbb{T}\cup\{(0,0)\}\cup\mathbb{T}^2=\mathbb{T}\times\overline{\mathbb{D}}\cup \overline{\mathbb{D}}\times\mathbb{T}\cup\{(0,0)\},$$
as in the statement.

{\bf Case G}: $j_0=-\infty , j'_0\in\mathbb{Z}, j_1=j_1'=\infty.$

This case is symmetric to the Case B.

{\bf Case H}: $j_0=-\infty, j_1=\infty , j'_0,j_1'\in\mathbb{Z}.$

This case is symmetric to the Case C.

{\bf Case I}: $j_0=j'_0=-\infty , j_1\in\mathbb{Z}, j'_1=\infty.$

This case is symmetric to Case E.

\end{proof}
\section{Taylor spectrum of the pair of isometries given by diagrams}
In this section we summarize previous results to obtain the main result of the paper.
\begin{theorem}\label{main}
  Let $M_w, M_z\in B(H_J)$, where $J$ is a non-simple diagram.

    Then $$\sigma_T(M_w,M_z)= \{(\mu,\lambda)\in \overline{\mathbb{D}}^2:
|\mu|^{\max(\rho_-,\rho_+)}\leq|\lambda|\leq |\mu|^{\min(\delta_-,\delta_+)}\},$$ 
where $\delta_-,\delta_+,\rho_-,\rho_+$ are defined by \eqref{deltarho+-} and inequalities with $0^0$ or $1^\infty$ are assumed to be satisfied.

\end{theorem}

\begin{proof}
Recall that isometries given by diagrams have empty point spectrum. Hence, by Remark \ref{T2iso}, we get $\Gamma_1(M_w,M_z)=\emptyset$ and in turn $\sigma_T(M_w,M_z)=\Gamma_2(M_w,M_z)\cup\Gamma_3(M_w,M_z).$ In particular, $\Gamma_2(M_w,M_z)\cup\Gamma_3(M_w,M_z)$ is closed.

Note that $\mathbb{T}^2$ lies in the closure of $\{(\mu,\lambda)\in \overline{\mathbb{D}}^2\setminus\mathbb{T}^2:
|\mu|^{\max(\rho_-,\rho_+)}\leq|\lambda|\leq |\mu|^{\min(\delta_-,\delta_+)}\}.$ Indeed, if $\max\{\rho_-,\rho_+\}=\infty,$ then $\mathbb{T}\times\mathbb{D}\subset \{(\mu,\lambda)\in \overline{\mathbb{D}}^2\setminus\mathbb{T}^2:
|\mu|^{\max(\rho_-,\rho_+)}\leq|\lambda|\leq |\mu|^{\min(\delta_-,\delta_+)}\}.$ If $\max\{\rho_-,\rho_+\}=0,$ then $\mathbb{D}\times\mathbb{T}\subset \{(\mu,\lambda)\in \overline{\mathbb{D}}^2\setminus\mathbb{T}^2:
|\mu|^{\max(\rho_-,\rho_+)}\leq|\lambda|\leq |\mu|^{\min(\delta_-,\delta_+)}\}.$ Otherwise, for any  $(\mu,\lambda)\in\mathbb{T}^2$ we may find a sequence $(\mu_n,\lambda_n)\in\mathbb{D}^2$ approximating $(\mu,\lambda),$ such that $|\lambda_n|=|\mu_n|^{\max\{\rho_-,\rho_+\}}$ for each $n$.

Summing up, the statement is equivalent to $$\Gamma_2(M_w,M_z)\cup\Gamma_3(M_w,M_z)\cup\mathbb{T}^2=\{(\mu,\lambda)\in \overline{\mathbb{D}}^2:
|\mu|^{\max(\rho_-,\rho_+)}\leq|\lambda|\leq |\mu|^{\min(\delta_-,\delta_+)}\},$$
where the set on the left is closed. 

Let us now check cases as in Theorem \ref{lem3cond}. Recall that $\delta_-\leq\eta_-\leq\rho_-$ and $\delta_+\leq\eta_+\leq\rho_+.$

If both $M_w, M_z$ have nontrivial Wold decomposition, then $\delta_-=\eta_-=\rho_-=0$ and $\delta_+=\eta_+=\rho_+=\infty.$ Hence $$\{(\mu,\lambda)\in \overline{\mathbb{D}}^2:
|\mu|^{\max(\rho_-,\rho_+)}\leq|\lambda|\leq |\mu|^{\min(\delta_-,\delta_+)}\}=\overline{\mathbb{D}}^2.$$

On the other hand, by Theorem \ref{lem2cond}  we get $$\Gamma_2(M_w,M_z)\setminus\mathbb{T}^2=\mathbb{D}^2$$
which by closeness of Taylor spectrum yields $$\sigma_T(M_w,M_z)=\overline{\mathbb{D}}^2.$$

Next, we consider the case where $M_w$ has nontrivial Wold decomposition and $M_z$ is a unilateral shift. Then $\delta_+=\eta_+=\rho_+=\infty$ and in turn
$$\{(\mu,\lambda)\in \overline{\mathbb{D}}^2:
|\mu|^{\max(\rho_-,\rho_+)}\leq|\lambda|\leq |\mu|^{\min(\delta_-,\delta_+)}\}=\{(\mu,\lambda)\in \overline{\mathbb{D}}^2:
|\lambda|\leq |\mu|^{\delta_-}\}.$$

By Theorems \ref{lem2cond} and \ref{lem3cond} we get
\begin{align*}\{(\mu,\lambda)\in (\mathbb{D}\setminus\{0\})^2 : & |\lambda|< |\mu|^{\eta_-} \}\cup(\mathbb{D}\times\{0\})\subset\Gamma_2(M_w,M_z)\setminus \mathbb{T}^2 \\&\subset\{(\mu,\lambda)\in (\mathbb{D}\setminus\{0\})^2 :  |\lambda|\leq |\mu|^{\eta_-} \}\cup\mathbb{D}\times\{0\}\end{align*}
and 
$$\Gamma_3(M_w,M_z)\cup\mathbb{T}^2=\{(\mu,\lambda)\in\overline{\mathbb{D}}^2: |\mu|^{\rho_-}\leq |\lambda|\leq |\mu|^{\delta_-} \textnormal{ or } |\mu|=1\}.$$

Since $\delta_-\leq\eta_-\leq \rho_-$ the inequalities $|\mu|^{\rho_-}\leq |\lambda|\leq |\mu|^{\delta_-}$ and $|\lambda|< |\mu|^{\eta_-}$ are equivalent to $|\lambda|\leq |\mu|^{\delta_-}$ on the joint domain $(\mathbb{D}\setminus\{0\})^2$ and we get $$\{(\mu,\lambda)\in (\mathbb{D}\setminus\{0\})^2:
|\lambda|\leq |\mu|^{\delta_-}\}\subset\Gamma_2(M_w,M_z)\cup \Gamma_3(M_w,M_z)\cup\mathbb{T}^2\subset \{(\mu,\lambda)\in \overline{\mathbb{D}}^2:
|\lambda|\leq |\mu|^{\delta_-}\}.$$ 
Since $\Gamma_2(M_w,M_z)\cup\Gamma_3(M_w,M_z)\cup\mathbb{T}^2$ is closed we get the statement.

The case where $M_w$ is a unilateral shift and $M_z$ has nontrivial Wold decomposition is symmetric.

For the last case, where  $M_w, M_z$ are unilateral shifts we get by Theorems \ref{lem2cond} and \ref{lem3cond} that
\begin{align*}\{(\mu,\lambda)\in (\mathbb{D}\setminus\{0\})^2 : &|\mu|^{\eta_+}< |\lambda|< |\mu|^{\eta_-} \}\cup\{(0,0)\}\subset\Gamma_2(M_w,M_z)\setminus \mathbb{T}^2 \\&\subset\{(\mu,\lambda)\in (\mathbb{D}\setminus\{0\})^2 : |\mu|^{\eta_+}\leq |\lambda|\leq |\mu|^{\eta_-} \}\cup\{(0,0)\}\end{align*}
and
$$\Gamma_3(M_w,M_z)\cup\mathbb{T}^2=\{(\mu,\lambda)\in\overline{\mathbb{D}}^2:
|\mu|^{\rho_-}\leq|\lambda|\leq |\mu|^{\delta_-}   \textnormal{ or }
|\mu|^{\delta_-}\leq|\lambda|\leq |\mu|^{\rho_+}  \textnormal{ or }
|\mu|^{\rho_+}\leq|\lambda|\leq |\mu|^{\delta_+}\}.$$

\begin{itemize}
    \item If $\delta_-\geq\rho_+,$ then $|\mu|^{\rho_-}\le|\mu|^{\eta_-}\le|\mu|^{\delta_-}\le|\mu|^{\rho_+}\le|\mu|^{\eta_+}\le|\mu|^{\delta_+}$ and hence $\Gamma_3(M_w,M_z)\cup\mathbb{T}^2=\{(\mu,\lambda)\in\overline{\mathbb{D}}^2:
|\mu|^{\rho_-}\leq|\lambda|\leq |\mu|^{\delta_+}\}$ contains $\Gamma_2(M_w,M_z)$ (which is thin as $|\mu|^{\eta_-}\leq|\mu|^{\eta_+}$, possibly reduced to $\{(0,0)\}$).
    \item If $\delta_+\geq\rho_-,$ then $|\mu|^{\rho_+}\le|\mu|^{\eta_+}\le|\mu|^{\delta_+}\le|\mu|^{\rho_-}\le|\mu|^{\eta_-}\le|\mu|^{\delta_-}$ and hence $\Gamma_2(M_w,M_z)\cup\Gamma_3(M_w,M_z)\cup\mathbb{T}^2=\{(\mu,\lambda)\in\overline{\mathbb{D}}^2:
|\mu|^{\rho_+}\leq|\lambda|\leq |\mu|^{\delta_-}\}.$
    \item If $\max\{\delta_-,\delta_+\}\leq\min\{\rho_-,\rho_+\},$ then $$\Gamma_3(M_w,M_z)\cup\mathbb{T}^2=\{(\mu,\lambda)\in\overline{\mathbb{D}}^2:
|\mu|^{\max\{\rho_-,\rho_+\}}\leq|\lambda|\leq |\mu|^{\min\{\delta_-,\delta_+\}}\}$$ and it contains $\Gamma_2(M_w, M_z)$.
\end{itemize}

\end{proof}

\section{Taylor spectrum of pairs of isometries - examples - continued}
In Remark \ref{ex_gp} we showed Taylor spectrum for pairs defined by the diagram being a subset of $\mathbb{Z}^2$ on the upper-right from some line. However technique used there works for quotient coefficient $-\frac{n}{m}$, as we need to define a unitary operator based on the integers $m,n$. Theorem \ref{main} extends the result to any line.  
\begin{remark}\label{oneSlope}
   If $\alpha>0,$ then $J=\{(i,j):j\geq-\alpha i\}$ is a diagram. Moreover, for such diagram $\rho_+=\rho_-=\delta_-=\delta_+=\frac1\alpha$. Therefore, $\sigma_T(M_w,M_z)=\{(\mu,\lambda)\in \overline{\mathbb{D}^2}:|\lambda|^\alpha=|\mu|\}. $  
\end{remark}
In the case above we get Taylor spectrum of Lebesgue measure $0$. On the other end we have diagrams $\mathbb{Z}_+^2$ and $\mathbb{Z}^2\setminus\mathbb{Z}_-^2$ where Taylor spectrum the whole $\cD^2.$ Such diagrams may be described as  subsets of $\mathbb{Z}^2$ on the upper-right from two perpendicular half-lines. The Theorem \ref{main} allows us to fill the gap between this two types of diagrams. More precisely, we consider diagrams defined as sets on the upper right from two half-lines. The size of Taylor spectrum depends on the coefficients of the half-lines.

\begin{remark}
    If $\beta>\alpha>0,$ then $$J=\{(i,j):j\geq-\alpha i, \textnormal{ for } i\geq 0 \textnormal{ and } j\geq-\beta i, \textnormal{ for } i< 0\}$$ is a diagram. Moreover, $\rho_+=\delta_+=\frac1\alpha$ and $\rho_-=\delta_-=\frac1\beta$. Therefore, $$\sigma_T(M_w,M_z)=\{(\mu,\lambda)\in\overline{\mathbb{D}^2}:|\lambda|^\beta\leq|\mu|\leq|\lambda|^\alpha\}. $$  
\end{remark}

\section*{Acknowledgement}

Research was supported by the Ministry of Science and Higher Education of the Republic
of Poland.

\bibliographystyle{alpha}

\end{document}